\documentclass[reqno,11pt]{amsart}

\usepackage{amsfonts,amsmath,amsthm,amssymb}
\usepackage{shortvrb,graphicx,psfrag}
\usepackage{wrapfig}
\usepackage[colorlinks,pdfpagelabels,pdfstartview = FitH,bookmarksopen
= true,bookmarksnumbered = true,linkcolor = blue,plainpages =
false,hypertexnames = false,citecolor = red,pagebackref=false]{hyperref}
\usepackage{relsize}
\usepackage{orcidlink}
\usepackage{xcolor}
\usepackage{appendix}
\usepackage{comment}
\usepackage{ulem}
\usepackage[makeroom]{cancel}
\usepackage[margin=1.4in]{geometry}
\usepackage{stmaryrd}
\usepackage{float}
\usepackage{mathrsfs}
\usepackage{mathabx}

\SetSymbolFont{stmry}{bold}{U}{stmry}{m}{n}

\allowdisplaybreaks
\usepackage{esint}

\allowdisplaybreaks[2]
\SetSymbolFont{stmry}{bold}{U}{stmry}{m}{n}
\let\TeXchi\chi
\newbox\chibox
\setbox0 \hbox{\mathsurround0pt $\TeXchi$}
\setbox\chibox \hbox{\raise\dp0 \box 0 }
\def\chi{\copy\chibox}

\DeclareMathOperator*{\esssup}{ess\,sup}

\DeclareMathOperator*{\essosc}{ess\,osc}
\DeclareMathOperator{\Tail}{Tail}
\DeclareMathOperator{\loc}{loc}

\newcommand{\R}{\mathbb{R}}

\newcommand{\dd}{\,\mathrm{d}}

\renewcommand{\rho}{\varrho}
\newcommand{\data}{\mathrm{data}}

\newcommand{\meas}[1]{\left|#1\right|}

\def\Xint#1{\mathchoice
    {\XXint\displaystyle\textstyle{#1}}
    {\XXint\textstyle\scriptstyle{#1}}
    {\XXint\scriptstyle\scriptscriptstyle{#1}}
    {\XXint\scriptscriptstyle\scriptscriptstyle{#1}}
    \!\int}
\def\XXint#1#2#3{\setbox0=\hbox{$#1{#2#3}{\int}$}
    \vcenter{\hbox{$#2#3$}}\kern-0.5\wd0}
\def\bint{\Xint-}
\def\dashint{\Xint{\raise4pt\hbox to7pt{\hrulefill}}}
\def\dashiint{\bint\kern-0.15cm\bint}

\newtheorem{theorem}{Theorem}[section]
\newtheorem{proposition}{Proposition}[section]
\newtheorem{lemma}{Lemma}[section]

\theoremstyle{definition}
\newtheorem{definition}{Definition}[section]
\newtheorem{remark}{Remark}[section]

\numberwithin{equation}{section}
\numberwithin{theorem}{section}
\numberwithin{proposition}{section}
\numberwithin{lemma}{section}
\numberwithin{remark}{section}
\title[Nonlocal nonlinear diffusion equations]{Continuity of solutions to\\ nonlocal nonlinear diffusion equations}

\author[Zh. Li]{Zheng Li\,\orcidlink{0009-0005-1696-1839}} 
\address{School of Mathematical Sciences\\ 
University of Electronic Science and Technology of China\\
2006 Xinyuan Avenue, 611731, Chengdu, P.R. China\\}           
\email{zhengli@uestc.edu.cn}                          
\thanks{}                                                

\date{}

\begin{document}



\begin{abstract}
We establish an explicit   modulus of continuity for locally bounded weak
solutions to nonlocal porous medium type equations with measurable
kernels. The nonlocal tail is only assumed to be integrable in time. The result lays a solid foundation for the intrinsic Harnack estimates of these equations.
\end{abstract}
\maketitle

\section{Introduction and Main Result}
Let $E\subset\R^N$ be open, $T>0$, $q>0$, and $s\in(0,1)$.
We consider locally bounded, local weak solutions to
\begin{equation}\label{equ01}
 \partial_t\boldsymbol{u}^q+\mathscr L u=0
 \qquad\text{in }E_T:=E\times(0,T),
\end{equation}
where $\boldsymbol{u}^q=|u|^{q-1}u$ for $u\ne0$ and
$\boldsymbol{0}^q=0$. The nonlocal operator is formally given by
\begin{equation*}
 \mathscr L u(x,t)=\operatorname{P.V.}
 \int_{\R^N}K(x,y,t)\bigl(u(x,t)-u(y,t)\bigr)\dd y .
\end{equation*}
The kernel is measurable, symmetric, and satisfies
\begin{equation}\label{kernel}
 \frac{C_o}{|x-y|^{N+2s}}\le K(x,y,t)=K(y,x,t)
 \le\frac{C_1}{|x-y|^{N+2s}}
\end{equation}
for almost every $(x,y,t)$, where $C_o,C_1>0$.
The parameters $\{N,s,q,C_o,C_1\}$ are referred to as the
$\data$.

The purpose of this paper is to give a quantitative local
continuity estimate when the nonlocal tail is merely
integrable in time. The dependence on its integrals over
shrinking time intervals is made explicit in the estimate
below.

\begin{theorem}[Local modulus of continuity]\label{main}
Let $u$ be a locally bounded, local weak solution to
\eqref{equ01}--\eqref{kernel} in $E_T$, in the sense of
Definition~\ref{def:weak}. There exist $C>1$ and
$\beta\in(0,s)$, depending only on the $\data$,
such that the following holds.
Let $Q=K_R(x_o)\times(T_1,T_2)\Subset E_T$, choose $M$
so that
\begin{equation}\label{main-amplitude}
 2M\ge\esssup_Q|u|+C\Tail(u;Q)^{1/q},
\end{equation}
and suppose
$(x_o,t_o)+Q_{2\rho}(\theta)\subset Q$ with $\theta=M^{q-1}$.
Then, for every $0<r<\rho$,
\begin{align}\label{main-modulus}
 &\essosc_{(x_o,t_o)+Q_r(\theta)}u\\ \notag
       &\quad\le CM\left(\frac r{\rho}\right)^{\beta/2}+CM\left[
       \frac1{M^q}\int_{t_o-\theta\rho^{2s}(r/\rho)^\beta}^{t_o}
       \int_{\R^N\setminus K_R(x_o)}
       \frac{|u(y,t)|}{|y-x_o|^{N+2s}}\dd y\dd t
                    \right]^{\min\{1,1/q\}}.
\end{align}
\end{theorem}

\begin{remark}
The nonlocal porous medium or fast diffusion equation is  often written as $\partial_t v + \mathscr L \boldsymbol{v}^m=0$. Thus $0<q<1$ corresponds to porous medium case,
whereas $q>1$ corresponds to fast diffusion. 
We allow solutions of either sign and require the equation
to hold only  locally  in $E_T$.
\end{remark}

\begin{remark}
Our result concerns from--$L_{\loc}^\infty$--to--$C_{\loc}^0$--type regularity. Whether solutions are locally bounded is a separate issue; see recent work \cite{DeFilippis2026} in this connection.
\end{remark}

\subsection{Novelty and Significance}

Nonlocal nonlinear diffusion equations arise, among other settings, in the boundary heat control problems formulated in the book \cite{Duvaut-Lions} of Duvaut and Lions. The continuity theory of their solutions was initiated by Athanasopoulos and Caffarelli \cite{Caff-10}. Over the past 15 years, they have been subject to intensive studies in various setups; see~\cite{ Bonforte-17, Bonf-Fig-Vaz-18, Bonf-Vaz-14, DPQRV-12,  maan-26} for instance. 
In particular, the regularity theory of these equations was further developed in
\cite{Caff-22,Caff-Vass-11, DPQR-16, DPQR-18, Kass-13};  they proved continuity of {\it globally bounded} weak
solutions.
 More recently, by combining techniques from \cite{Liao-21, Liao-24a},
 H\"older regularity of {\it locally bounded} solutions with the $L^\infty$-in-time tail is proven in \cite{KLP-26}, whereas Q.~Li \cite{Li-26} established it under the weaker $L^{1+\varepsilon}$-in-time tail. 

The present work establishes local continuity at the natural $L^1$-in-time tail endpoint. More specifically, we consider locally bounded, possibly sign-changing weak solutions, require the equation to hold locally, and allow measurable symmetric kernels depending on both space and time. The result applies to every $q>0$, covering both porous medium and fast diffusion cases. No global boundedness or additional time integrability of the tail is imposed. 

The general modulus of continuity \eqref{main-modulus} is a natural continuation of the one derived in \cite{Liao-24b} for the parabolic fractional $p$-Laplacian.
It separates a power decay arising from the local oscillation argument from the contribution of the exterior values over a shrinking time interval. The latter vanishes by absolute continuity of the tail integral, without any prescribed rate. Thus the estimate describes how interior continuity depends on the temporal concentration of the solution outside the reference domain. When the tail belongs to $L^{1+\varepsilon}$ in time, the same estimate yields a Hölder modulus. This kind of Hölder regularity generalizes the linear result in \cite[Theorem~1.5]{Kass-23}; see also \cite{Byun-ampa,Liao-24b} for the fractional parabolic $p$-Laplacian.

A central issue is the interaction between the control of the exterior values and the degeneracy or singularity of the time nonlinearity at zero. The intrinsic time scale depends on the amplitude of the solution, and hence changes during the oscillation reduction. The tail must therefore be controlled on time intervals whose lengths are determined by the same iteration. Theorem \ref{main} shows that these two effects can be treated together, with constants depending only on the data and with the exterior contribution displayed explicitly. We will come back to the technical aspect later.

The continuity conclusion also supports the pointwise theory of these equations. In particular, it provides the continuous representative and the topological information required in the intrinsic Harnack argument of \cite[Theorem~7.1]{Liao-26}; see also~\cite[Remark 1.3]{Liao-26}. This connection makes endpoint continuity relevant to the development of intrinsic Harnack estimates under natural tail assumptions.

\subsection{Method of proof}

As in the local theory of porous medium type equations, we distinguish cylinders where the solution approaches or crosses zero from cylinders where it is quantitatively separated from zero; cf.~\cite{DB-Fried-85, DBGV-mono, Liao-21}. The distinction is necessary because adding constants to a solution does not yield another one due to the power nonlinearity in the time derivative. To distinguish these cylinders, the measure concentration of $u$ near zero is at stake. If it always has a sizable portion near zero along nested, intrinsic cylinders, then we keep reducing its amplitude. Otherwise, we show that $u$ must be away from zero in the next nested cylinder, where then the non-degenerate regularity theory can be applied after normalization. 

Apart from this main idea, the treatment of the nonlocal tail in this argument is quite delicate. It enters into positivity estimates in the following pattern: either the tail is large relative to the current amplitude, or the solution satisfies an improved pointwise bound in a smaller cylinder. 

With the positivity tools at hand, we then construct nested cylinders and an accompanying sequence of amplitude bounds. The amplitudes are allowed to decrease only as far as permitted by the tail. The exterior integral at each step is split into the region outside the original reference domain and the intervening spatial annuli. On the annuli, bounds obtained at earlier steps produce a geometric sum that can be absorbed into the amplitude reduction. The remaining contribution is the original exterior tail integrated over the current shrinking time interval. Preserving this interval throughout the iteration is essential: replacing it simply by the initial time interval would leave an error that need not tend to zero, cf.~\eqref{discrete1}.

If the solution becomes separated from zero at some stage, we stop the near-zero iteration and normalize at that scale. The normalized solution lies in a fixed positive range. The endpoint oscillation argument of \cite{Liao-24b} can then be adapted to this normalized equation.  

The final step is to express both alternatives in the original intrinsic geometry and obtain an estimate independent of the stopping index. For a prescribed reference radius, we compare the stopping scale with an intermediate iteration scale. If the iteration reaches that intermediate scale, its amplitude estimate controls the oscillation. If it stops earlier, the normalized nondegenerate estimate applies. Tracking the amplitude factors, the annular contributions, and the time interval in the rescaled tail places both estimates on a common cylinder. This quantitative matching yields the explicit modulus in \eqref{main-modulus}.

\subsection{Definitions and notation}
For $x_o\in\R^N$, let $K_\rho(x_o)$ denote the ball of
radius $\rho$ centered at $x_o$, and write $K_\rho=K_\rho(0)$.
We use the cylinders
\begin{align*}
 (x_o,t_o)+Q_\rho(\theta)&=K_\rho(x_o)\times(t_o-\theta\rho^{2s},t_o),\\
  (x_o,t_o)+Q(R,S)&=K_R(x_o)\times(t_o-S,t_o).
\end{align*}
If the context is clear, we omit $(x_o,t_o)$.
For $Q=K_R(x_o)\times I$, set
\begin{equation*}
 \Tail(f;Q)=\int_I\int_{\R^N\setminus K_R(x_o)}
       \frac{|f(y,t)|}{|y-x_o|^{N+2s}}\dd y\dd t .
\end{equation*}

Unless otherwise specified,
$C\ge1$ depends only on
the $\data$ and may change from line to line.
The constants $c_q,C_q>0$ depend only on $q$.
For a measurable set $A$, $|A|$ denotes its Lebesgue measure.\\
For simplicity, with $a\in \mathbb{R}$ and $\alpha>0$, we define
\[
\boldsymbol{a}^\alpha=|a|^{\alpha-1}a \text{ for } a\neq 0 \,\,\text{and}\,\,\boldsymbol{0}^\alpha=0. 
\]
We use the following notion of weak solution.
\begin{definition}\label{def:weak}
A measurable function $u:\R^N\times(0,T)\to\R$ is a local weak
subsolution to \eqref{equ01} if
\[
 u\in C_{\loc}\bigl((0,T);L^{q+1}_{\loc}(E)\bigr)
 \cap L^2_{\loc}\bigl((0,T);W^{s,2}_{\loc}(E)\bigr),
\]
if, for every compact interval $[t_1,t_2]\subset(0,T)$,
\begin{equation*}
 \int_{t_1}^{t_2}\int_{\R^N}
 \frac{|u(x,t)|}{(1+|x|)^{N+2s}}\dd x\dd t<\infty,
\end{equation*}
and if
\begin{align}\label{weak}
 &\left.\int_D\boldsymbol{u}^q\phi\dd x\right|_{t_1}^{t_2}
 -\int_{t_1}^{t_2}\int_D
       \boldsymbol{u}^q\partial_t\phi\dd x\dd t\notag\\
 &\quad+\frac12\int_{t_1}^{t_2}\int_{\R^N}\int_{\R^N}
 (u(x,t)-u(y,t))(\phi(x,t)-\phi(y,t))
 K(x,y,t)\dd y\dd x\dd t\le0
\end{align}
for nonnegative admissible test functions with compact spatial
support in $D\Subset E$. Here the test functions belong to
\[
 L^2\bigl(t_1,t_2;W^{s,2}_0(D)\bigr)
 \cap W^{1,q+1}\bigl(t_1,t_2;L^{q+1}(D)\bigr),
\]
and are extended by zero outside $D$.
A supersolution satisfies the reverse inequality. A weak solution
is both a subsolution and a supersolution.
\end{definition}

\subsection{Structure of the Paper}
Section~\ref{sec:prelim} collects the preliminary inequalities.
The energy estimate with time-dependent truncations is proved in
Section~\ref{sec:energy}. Sections~\ref{sec:away-dg} and
\ref{sec:near} contain the measure estimates used in the
intrinsic iteration. The proof of Theorem~\ref{main} is given
in Section~\ref{sec:continuity}.

\

\noindent\textbf{Acknowledgements.} The author would like to express his sincere gratitude to Professor Naian Liao, for suggesting the problem and for his generous guidance throughout this work.

\section{Preliminaries}\label{sec:prelim}
We collect the inequalities used below; see
\cite{Liao-24a,Misawa-25}. 

\subsection{Algebraic inequalities}
The signed-power estimate below is
\cite[Lemma 2.9]{Misawa-25}. We also record the corresponding
estimates for the time-energy primitive.
\begin{lemma}
For $q>0$ and $a,b\in\R$,
\begin{equation}\label{power-difference}
 c_q(|a|+|b|)^{q-1}|a-b|
 \le|\boldsymbol{a}^q-\boldsymbol{b}^q|
 \le C_q(|a|+|b|)^{q-1}|a-b|.
\end{equation}
At $a=b=0$, all expressions are assigned the value zero.
For
\begin{equation*}
 g_\pm(a,k)=\pm\int_k^a |r|^{q-1}(r-k)_\pm\dd r,
\end{equation*}
one has $g_\pm\ge0$ and
\begin{equation}\label{g-comparison}
 c_q(|a|+|k|)^{q-1}(a-k)_\pm^2
 \le g_\pm(a,k)
 \le C_q(|a|+|k|)^{q-1}(a-k)_\pm^2.
\end{equation}
\end{lemma}

In particular, if $|a|\le2M$ and $M/4\le k\le M$, then
\begin{align}
 c_qM^{q-1}(a-k)_-^2
 &\le g_-(a,k)\le C_qM^{q-1}(a-k)_-^2,
       \label{g-lower-level}\\
 (\boldsymbol{a}^q-\boldsymbol{k}^q)_-
 &\ge c_qM^{q-1}(a-k)_-.\label{power-lower-level}
\end{align}
For $M\le k\le2M$, the analogous estimates hold for
$g_+(a,k)$ and $(\boldsymbol{a}^q-\boldsymbol{k}^q)_+$.
We shall also use
\begin{equation*}
 qg_-(\boldsymbol{z}^{1/q},k)
 =\int_z^{\boldsymbol{k}^q}
       (k-\boldsymbol{r}^{1/q})_+\dd r
\end{equation*}
and
\begin{align}\label{g-continuity}
 q\left|g_\pm(\boldsymbol{z_1}^{1/q},k)
          -g_\pm(\boldsymbol{z_2}^{1/q},k)\right|\le
 (|z_1|^{1/q}+|z_2|^{1/q}+|k|)|z_1-z_2|.
\end{align}

\subsection{Fractional embedding and auxiliary estimates}
Write
\[
 [v]_{W^{s,2}(K_\rho)}^2
 =\int_{K_\rho}\int_{K_\rho}
       \frac{|v(x)-v(y)|^2}{|x-y|^{N+2s}}\dd y\dd x .
\]
The next statement is the case $p=2$ of \cite[Proposition A.3]{Liao-24a}; see also \cite[Lemma 2.2]{Misawa-25}.

\begin{lemma}\label{embedding}
Let $I$ be a bounded interval and
$v\in L^\infty(I;L^2(K_\rho))\cap L^2(I;W^{s,2}(K_\rho))$, which is compactly supported in $K_{(1-\sigma)\rho}$ for $\sigma\in (0,1)$ and for a.e. $t\in I$.  Then

\begin{align*}
 \int_I\int_{K_\rho}|v|^{\widehat m}\dd x\dd t
 &\le C(N,s)
 \left(\esssup_{t\in I}\bint_{K_\rho}|v(x,t)|^2\dd x\right)^{2s/N}
 \notag\\
 &\quad\times
 \int_I\left(\rho^{2s}[v(\cdot,t)]_{W^{s,2}(K_\rho)}^2
             +\sigma^{-(N+2s)}\int_{K_\rho}|v|^2\dd x\right)\dd t,
\end{align*}
where $\widehat m=2(1+2s/N)$.
\end{lemma}
Now, we collect some useful auxiliary estimates to be used.
Let $\zeta\in C^1(K_\rho,[0,1])$,
 and $|D\zeta|\le L$, then for $x\in K_\rho$,
\begin{equation}\label{cutoff-integral}
 \int_{K_\rho}\frac{|\zeta(x)-\zeta(y)|^2}{|x-y|^{N+2s}}\dd y
 \le C(N,s)L^2\rho^{2-2s}.
\end{equation}
Moreover,
\begin{align}\label{product-cutoff}
 |v(x)\zeta(x)-v(y)\zeta(y)|^2
 &\le2\zeta(x)^2|v(x)-v(y)|^2
 +2v(y)^2|\zeta(x)-\zeta(y)|^2.
\end{align}
For $x\in K_\rho(x_o)$, $y\notin K_R(x_o)$, and $R\ge2\rho$,
\begin{equation}\label{tail-distance}
 |x-y|\ge\tfrac12|y-x_o|,\qquad
 \int_{\R^N\setminus K_\rho(x_o)}
       \frac{\dd y}{|y-x_o|^{N+2s}}=\frac{C(N,s)}{\rho^{2s}}.
\end{equation}
More generally, if $x\in K_{\hat\rho}(x_o)$,
$y\notin K_\rho(x_o)$, and $\hat\rho<\rho$, then
\begin{equation}\label{separation}
 \frac{|y-x_o|}{|x-y|}\le\frac{\rho}{\rho-\hat\rho}.
\end{equation}

\subsection{Time regularization and geometric convergence}
We use the forward and backward exponential regularizations
as in \cite[Appendix B]{Liao-24a} and
\cite[Lemma B.1]{Misawa-25}. For a Banach space $X$ and $a<b$,
set
\begin{equation}\label{mollifiers}
 (f)_h(t)=\frac1h\int_a^t e^{(\tau-t)/h}f(\tau)\dd\tau,\qquad
 [f]_h(t)=\frac1h\int_t^b e^{(t-\tau)/h}f(\tau)\dd\tau .
\end{equation}
Both regularizations converge strongly to $f$ in
$L^p(a,b;X)$ for $1\le p<\infty$, and
\begin{equation}\label{mollifier-derivatives}
 \partial_t(f)_h=\frac{f-(f)_h}{h},\qquad
 \partial_t[f]_h=\frac{[f]_h-f}{h}.
\end{equation}
If $f\in C([a,b];X)$ and $f(b)=0$, then by \cite[Lemma B.1, (vi)]{Misawa-25}
\begin{equation}\label{backward-uniform}
 [f]_h\longrightarrow f\quad\text{in }C([a,b];X).
\end{equation}
The condition $f(b)=0$ is essential here, since $[f]_h(b)=0$.

We finally recall the fast geometric convergence lemma;
see \cite[Lemma 2.10]{Misawa-25}.
\begin{lemma}\label{fast}
Let $Y_n\ge0$ satisfy $Y_{n+1}\le AB^nY_n^{1+\alpha}$,
where $A>0$, $B>1$, and $\alpha>0$.
If
\begin{equation*}
 Y_0\le A^{-1/\alpha}B^{-1/\alpha^2},
\end{equation*}
then $Y_n\to0$ as $n\to\infty$.
\end{lemma}

\section{The energy estimate}\label{sec:energy}
\begin{proposition}\label{energy}
Let $u$ be a local weak subsolution (respectively, supersolution)
of \eqref{equ01}, and let $Q(R,S)=K_R\times(t_o-S,t_o)\Subset E_T$.
Let $k\in W^{1,1}(t_o-S,t_o)$, and let $0\le\zeta\le1$ be a
piecewise smooth cutoff whose spatial support is contained in
a fixed compact subset of $K_R$. Write $w_\pm=(u-k(t))_\pm$.
With the upper sign for a subsolution and the lower sign for a
supersolution, the following estimate holds:
\begin{align}\label{energy-full}
 &\left.q\int_{K_R}g_\pm(u,k)\zeta^2\dd x\right|_{t_o-S}^{t_o}
 \ \pm\int_{t_o-S}^{t_o}\int_{K_R}k'(t)\zeta^2
             (\boldsymbol{u}^q-\boldsymbol{k}^q)_\pm\dd x\dd t
 \notag\\
 &\quad+\frac12\int_{t_o-S}^{t_o}\int_{K_R}\int_{K_R}
      |\zeta(x,t)w_\pm(x,t)-\zeta(y,t)w_\pm(y,t)|^2K(x,y,t)\dd y\dd x\dd t
 \notag\\
 &\quad+\int_{t_o-S}^{t_o}\int_{K_R}\int_{K_R}
      w_\mp(y,t)w_\pm(x,t)
      \zeta(x,t)^2 K(x,y,t)\dd y\dd x\dd t
 \notag\\
 &\le q\int_{t_o-S}^{t_o}\int_{K_R}g_\pm(u,k)\partial_t(\zeta^2)\dd x\dd t
 \notag\\
 &\quad+\int_{t_o-S}^{t_o}\int_{K_R}\int_{K_R}
      |\zeta(x,t)-\zeta(y,t)|^2
      w_\pm(x,t)\,w_\pm(y,t)
      K(x,y,t)\dd y\dd x\dd t
 \notag\\
 &\quad+\int_{t_o-S}^{t_o}\int_{K_R}\int_{\R^N\setminus K_R}
      w_\pm(y,t)\zeta(x,t)^2w_\pm(x,t)
      K(x,y,t)\dd y\dd x\dd t .
\end{align}
The same statement holds with any intermediate terminal time
in place of $t_o$.
\end{proposition}

\begin{proof}
We prove the supersolution case.

\smallskip\noindent\textit{Step 1: the test function and the time term.}
Choose $0<\varepsilon<S/4$ and set
\[
 \psi_\varepsilon(t)=
 \begin{cases}
 (t-t_o+S)/\varepsilon,&t_o-S\le t<t_o-S+\varepsilon,\\
 1,&t_o-S+\varepsilon\le t\le t_o-\varepsilon,\\
 2(t_o-\varepsilon/2-t)/\varepsilon,
       &t_o-\varepsilon<t<t_o-\varepsilon/2,\\
 0,&t_o-\varepsilon/2\le t\le t_o.
 \end{cases}
\]
Set $\varphi=\zeta^2\psi_\varepsilon w_-$ and use
$\phi=[\varphi]_h$ in \eqref{weak}, with the regularizations
in \eqref{mollifiers} taken over $(t_o-S,t_o)$.
Fubini's theorem gives 
\begin{align*}
  -\int_{t_o-S}^{t_o}\int_{K_R}
       \boldsymbol{u}^q\partial_t[\varphi]_h\dd x\dd t =
 \int_{t_o-S}^{t_o}\int_{K_R}
       \varphi\,\partial_t(\boldsymbol{u}^q)_h\dd x\dd t .
\end{align*}
Here $[\varphi]_h(t_o)=0$. The initial endpoint contribution is
\[
 -\int_{K_R}\boldsymbol{ u}^q(x,t_o-S)
                 [\varphi]_h(x,t_o-S)\dd x.
\]
It tends to zero as
$h\downarrow0$, by \eqref{backward-uniform},
H\"older's inequality, and $\varphi(t_o-S)=0$.

The map $z\mapsto(k-\boldsymbol{z}^{1/q})_+$ is nonincreasing.
Thus 
\[
 \left[(u-k)_-
       -\boldsymbol{\bigl((\boldsymbol{u}^q)_h}^{1/q}-k\bigr)_-
 \right]\left[\boldsymbol{u}^q-(\boldsymbol{u}^q)_h\right]
 \le0.
\]
Using \eqref{mollifier-derivatives}, we infer
\begin{align}\label{time-monotonicity}
 \int_{t_o-S}^{t_o}\int_{K_R}
       \varphi\,\partial_t(\boldsymbol{u}^q)_h\dd x\dd t\le
 \int_{t_o-S}^{t_o}\int_{K_R}
       \zeta^2\psi_\varepsilon
       \bigl(k-\boldsymbol{(\boldsymbol{u}^q)_h}^{1/q}\bigr)_+
       \partial_t(\boldsymbol{u}^q)_h\dd x\dd t.
\end{align}
By the formal differentiation, we have
\begin{align*}
 \bigl(k-\boldsymbol{(\boldsymbol{u}^q)_h}^{1/q}\bigr)_+
             \partial_t(\boldsymbol{u}^q)_h=-q\,\partial_t
       g_-\bigl(\boldsymbol{(\boldsymbol{u}^q)_h}^{1/q},k\bigr)
       +k'(t)\bigl(\boldsymbol{k}^q-(\boldsymbol{u}^q)_h\bigr)_+ .
\end{align*}
Substituting it into \eqref{time-monotonicity}, integrating by parts, using
$\psi_\varepsilon(t_o-S)=\psi_\varepsilon(t_o)=0$,
turns the right-hand side of \eqref{time-monotonicity} into
\begin{align}\label{time-three}
 &q\int_{t_o-S}^{t_o}\int_{K_R}
       g_-\bigl(\boldsymbol{(\boldsymbol{u}^q)_h}^{1/q},k\bigr)
       \zeta^2\psi_\varepsilon'\dd x\dd t\notag\\
 &\quad+q\int_{t_o-S}^{t_o}\int_{K_R}
       g_-\bigl(\boldsymbol{(\boldsymbol{u}^q)_h}^{1/q},k\bigr)
       \psi_\varepsilon\partial_t(\zeta^2)\dd x\dd t\notag\\
 &\quad+\int_{t_o-S}^{t_o}\int_{K_R}
       \zeta^2\psi_\varepsilon k'(t)
       \bigl(\boldsymbol{k}^q-(\boldsymbol{u}^q)_h\bigr)_+
       \dd x\dd t .
\end{align}

\smallskip\noindent\textit{Step 2: passage to the limit in time.}
We first send $h\downarrow0$ in \eqref{time-three}. Since $u\in C([t_o-S,t_o];L^{q+1}(K_R))$,
the signed power $\boldsymbol{u}^q$ belongs to
$C([t_o-S,t_o];L^{(q+1)/q}(K_R))$.
Its forward averages converge strongly in
$L^{(q+1)/q}(K_R\times(t_o-S,t_o))$ and in
$L^{(q+1)/q}(K_R)$ at each $t>t_o-S$.
They are uniformly bounded in
$L^\infty(t_o-S,t_o;L^{(q+1)/q}(K_R))$.
By \eqref{g-continuity} and H\"older's inequality,
\begin{equation}\label{g-limit}
 g_-\bigl(\boldsymbol{(\boldsymbol{u}^q)_h}^{1/q},k\bigr)
       \longrightarrow g_-(u,k)
       \quad\text{in }L^1(K_R\times(t_o-S,t_o)).
\end{equation}
Here $k$ is bounded and
$\boldsymbol{(\boldsymbol{u}^q)_h}^{1/q}$ is uniformly
bounded in $L^{q+1}$ on this cylinder.
For fixed $\varepsilon$, \eqref{g-limit} gives the limits
of the first two terms in \eqref{time-three}.

For the last term, the absolute value of the error is at most
\[
 \int_{t_o-S}^{t_o}|k'(t)|
       \|(\boldsymbol{u}^q)_h(t)-\boldsymbol{u}^q(t)\|
                     _{L^1(K_R)}\dd t .
\]
The integrand converges to zero for $t>t_o-S$ and is bounded
by a constant times $|k'(t)|$. The dominated convergence theorem
forces the above integral to become zero. Hence, after $h\downarrow0$, we get from \eqref{time-three} that
\begin{align*}
 &q\int_{t_o-S}^{t_o}\int_{K_R}
       g_-\bigl(u,k\bigr)
       \zeta^2\psi_\varepsilon'\dd x\dd t
 +q\int_{t_o-S}^{t_o}\int_{K_R}
       g_-\bigl(u,k\bigr)
       \psi_\varepsilon\partial_t(\zeta^2)\dd x\dd t\\
&\quad+\int_{t_o-S}^{t_o}\int_{K_R}
       \zeta^2\psi_\varepsilon k'(t)
       \bigl(\boldsymbol{k}^q-\boldsymbol{u}^q\bigr)_+
       \dd x\dd t .
\end{align*}

Next let $\varepsilon\downarrow0$. By \eqref{g-continuity},
the function $t\mapsto\int_{K_R}g_-(u,k)\zeta^2\dd x$ is
continuous. The limit of \eqref{time-three} is therefore
\begin{align}\label{time-limit}
 &-q\left.\int_{K_R}g_-(u,k)\zeta^2\dd x\right|_{t_o-S}^{t_o}
      +q\int_{t_o-S}^{t_o}\int_{K_R}
       g_-(u,k)\partial_t(\zeta^2)\dd x\dd t\notag\\
 &\quad+\int_{t_o-S}^{t_o}\int_{K_R}
       k'(t)\zeta^2(\boldsymbol{k}^q-\boldsymbol{u}^q)_+
       \dd x\dd t .
\end{align}

\smallskip\noindent\textit{Step 3: convergence of the spatial term.}
We first split the fractional term into

\begin{align}\label{spatial_term}
&\int_{t_{o}-S}^{t_{o}} \int_{K_{R}} \int_{K_{R}}(u(x, t)-u(y, t))([\varphi]_h(x, t)-[\varphi]_h(y, t)) K(x, y, t) \dd x \dd y \dd t\notag\\
&+2 \int_{t_{o}-S}^{t_{o}} \int_{\mathbb{R}^{N}\setminus K_{R}} \int_{K_{R}}u(x, t)[\varphi]_h(x, t) K(x, y, t) \dd x \dd y \dd t\notag\\
&-2 \int_{t_{o}-S}^{t_{o}} \int_{\mathbb{R}^{N}\setminus K_{R}} \int_{K_{R}}u(y, t)[\varphi]_h(x, t) K(x, y, t) \dd x \dd y \dd t.
\end{align}

On $K_R\times K_R$, the Cauchy--Schwarz inequality bounds the
error by
\begin{align*}
 C_1 
 \left(\int_{t_o-S}^{t_o} [u(\cdot,t)]_{W^{s,2}(K_R)}^2\dd t\right)^{1/2} 
 \left(\int_{t_o-S}^{t_o}[[\varphi]_h(\cdot,t)-\varphi(\cdot,t)]
                         _{W^{s,2}(K_R)}^2\dd t\right)^{1/2}
 \stackrel{h\downarrow 0}{\longrightarrow}0.
\end{align*}
Thus the first term of \eqref{spatial_term} goes to
\[
\int_{t_{o}-S}^{t_{o}} \int_{K_{R}} \int_{K_{R}}(u(x, t)-u(y, t))(\varphi(x, t)-\varphi(y, t)) K(x, y, t) d x d y d t.
\]
For the interactions with $\R^N\setminus K_R$, use the positive
distance between the fixed spatial support of $\zeta$ and
$\partial K_R$.
The error of the second term of \eqref{spatial_term} tends to zero by local $L^2$
convergence. The error of the final term is bounded by
\begin{align*}
 C(\zeta,R)&\,
 \sup_{t\in[t_o-S,t_o]}
   \|[\varphi]_h(t)-\varphi(t)\|_{L^1(K_R)} \int_{t_o-S}^{t_o}\int_{\R^N}
          \frac{|u(y,t)|}{(1+|y|)^{N+2s}}\dd y\dd t
 \longrightarrow0 .
\end{align*}
The supremum tends to zero by \eqref{backward-uniform};
$\varphi$ vanishes near $t_o$.
Finally $\psi_\varepsilon\to1$ almost everywhere.
Thus, \eqref{spatial_term} tends to
\[
\int_{t_{o}-S}^{t_{o}} \int_{\mathbb{R}^{N}} \int_{\mathbb{R}^{N}}(u(x, t)-u(y, t))\left(\zeta^{2} w_{-}(x, t)-\zeta^{2} w_{-}(y, t)\right) K(x, y, t) \dd x \dd y \dd t ,
\]
which for simplicity we denote as $\mathbf{I}$.
The next step consists of using an algebraic identity to decompose the above integral.

\smallskip\noindent\textit{Step 4: the truncation algebra.}
Fix $t$ and use $w_\pm=(u-k)_\pm$.
The exact identity is
\begin{align*}
 (u(x)&-u(y))\bigl(\zeta(x)^2w_-(x)-\zeta(y)^2w_-(y)\bigr)\notag\\
 =&-|\zeta(x)w_-(x)-\zeta(y)w_-(y)|^2
      +w_-(x)w_-(y)|\zeta(x)-\zeta(y)|^2\notag\\
 &-w_+(x)w_-(y)\zeta(y)^2-w_+(y)w_-(x)\zeta(x)^2 .
\end{align*}
Then, multiplying above identity by $\frac12K(x,y,t)$,
 integrating it over $(t_o-S,t_o)\times\mathbb{R}^N\times \mathbb{R}^N$, and then decomposing the integrals yield
\begin{align}\label{five_spacial}
\frac{\mathbf{I}}{2}
=&-\frac{1}{2}\int_{t_{o}-S}^{t_{o}} \int_{K_R} \int_{K_R}|\zeta(x,t)w_-(x,t)-\zeta(y,t)w_-(y,t)|^2 K(x, y, t) \dd x \dd y \dd t \notag\\
&-\int_{t_{o}-S}^{t_{o}} \int_{\mathbb{R}^{N}\setminus K_R} \int_{\mathbb{R}^{N}}|\zeta(x,t)w_-(x,t)|^2 K(x, y, t) \dd x \dd y \dd t\notag \\
&+\frac{1}{2}\int_{t_{o}-S}^{t_{o}} \int_{K_R} \int_{K_R}w_-(x,t)w_-(y,t)|\zeta(x,t)-\zeta(y,t)|^2 K(x, y, t) \dd x \dd y \dd t\notag\\
&+\int_{t_{o}-S}^{t_{o}} \int_{\mathbb{R}^{N}\setminus K_R} \int_{\mathbb{R}^{N}}w_-(x,t)w_-(y,t)\zeta(x,t)^2 K(x, y, t) \dd x \dd y \dd t\notag\\
&-\int_{t_{o}-S}^{t_{o}} \int_{K_R} \int_{K_R}w_+(x,t)w_-(y,t)\zeta(y,t)^2 K(x, y, t) \dd x \dd y \dd t\notag\\
&-\int_{t_{o}-S}^{t_{o}} \int_{K_R} \int_{ \mathbb{R}^N\setminus K_R}w_+(x,t)w_-(y,t)\zeta(y,t)^2 K(x, y, t) \dd x \dd y \dd t.
\end{align}

Dropping the second and the last terms of \eqref{five_spacial} and combining the remaining terms with \eqref{time-limit}, we prove
\eqref{energy-full}. The subsolution case is obtained by testing
with $(u-k)_+$, or by applying the supersolution result to $-u$.
\end{proof}

\begin{remark}\label{energy-use}
If $\zeta$ vanishes at the initial time, apply
\eqref{energy-full} up to an arbitrary terminal time, and then
take the essential supremum. Adding the resulting bound to
the bound for the full time interval controls both the
time supremum and the spatial energy, at the cost of a fixed
constant. When a moving level is used to absorb the 
tail, this absorption is performed before taking the
supremum. The coefficient $1/2$ is harmless in all subsequent
iterations, but it must agree with the normalization of
\eqref{weak}.
\end{remark}

\section{A De Giorgi lemma separating the solution from zero}
\label{sec:away-dg}
In this and the next sections, fix an reference cylinder
\[
 Q=K_R(x_o)\times(T_1,T_2)\Subset E_T.
\]
Throughout the sections, we assume for some $M>0$,
\[
 \esssup_{Q}|u|\le2M 
\]
and denote $\theta=M^{q-1}$.
Every time interval used in a lemma is assumed to lie in
$(T_1,T_2)$. This convention is necessary also for the
forward time intervals in Section~\ref{sec:near}.

\begin{lemma}\label{De01}
There exist $\nu\in(0,1)$ and $\gamma_{**}\ge1$, depending only on
$\data$, with the following property. Let $u$ be a locally
bounded, local weak supersolution, let $2\rho\le R$, and assume
$(x_o,t_o)+Q_\rho(\theta)\subset Q$.
If
\[
 |\{u\le M\}\cap((x_o,t_o)+Q_\rho(\theta))|
 \le\nu|Q_\rho(\theta)|,
\]
then either
\begin{equation}\label{dg-tail}
 \gamma_{**}\Tail(u_-;Q)>M^q
\end{equation}
or
\begin{equation}\label{dg-conclusion}
 u\ge M/4\qquad\text{a.e. in }
 (x_o,t_o)+Q_{\rho/2}(\theta).
\end{equation}
For a subsolution, the corresponding statement applies to $-u$,
with $u_+$ in \eqref{dg-tail}.
\end{lemma}
\begin{proof}
Translate $(x_o,t_o)$ to $(0,0)$.
We assume that \eqref{dg-tail} fails and choose $\gamma_{**}$
shortly. Introduce the levels and radii
\begin{align*}
 k_n&=\frac M2+\frac{M}{2^{n+1}},&
 \rho_n&=\frac\rho2+\frac{\rho}{2^{n+1}},\\
 \hat\rho_n&=\frac{3\rho_n+\rho_{n+1}}4,&
 \tilde\rho_n&=\frac{\rho_n+\rho_{n+1}}2,&
 \bar\rho_n&=\frac{\rho_n+3\rho_{n+1}}4 .
\end{align*}
Set $$
\left\{\begin{array}{l}
K_{n}=K_{\rho_{n}}, \quad \widetilde{K}_{n}=K_{\tilde{\rho}_{n}}, \quad \widehat{K}_{n}=K_{\hat{\rho}_{n}}, \quad \overline{K}_{n}=K_{\bar{\rho}_{n}}, \\[2pt]
Q_{n}=K_n \times\left(-\theta \rho_{n}^{2 s}, 0\right), \quad \widetilde{Q}_{n}=\widetilde{K}_{n} \times\left(-\theta \tilde{\rho}_{n}^{2 s}, 0\right), \\[2pt]
\widehat{Q}_{n}=\widehat{K}_{n} \times\left(-\theta \hat{\rho}_{n}^{2 s}, 0\right), \quad \overline{Q}_{n}=\overline{K}_{n} \times\left(-\theta \bar{\rho}_{n}^{2 s}, 0\right).
\end{array}\right.
$$
Thus
\[
 Q_{n+1}\subset\overline Q_n\subset\widetilde Q_n
                 \subset\widehat Q_n\subset Q_n.
\]
Choose $\zeta$ equal to one in $\widetilde Q_n$, supported
spatially in $\widehat K_n$, vanishing at the bottom of $Q_n$,
and zero at times below the bottom of $\widehat Q_n$.
It can be chosen so that
\begin{equation*}
 |D\zeta|\le\frac{C2^n}{\rho},\qquad
 |\partial_t\zeta|\le\frac{C(s)2^n}{\theta\rho^{2s}}.
\end{equation*}

\smallskip\noindent\textit{Step 1: moving levels and tail cancellation.}
Define
\begin{equation*}
 \ell(t)=\widetilde\gamma M^{1-q}\int_{T_1}^t
       \int_{\R^N\setminus K_R}
       \frac{u_-(y,\tau)}{|y|^{N+2s}}\dd y\dd\tau,\qquad
 \mathsf{k}_n(t)=k_n-\ell(t).
\end{equation*}
Here $\widetilde\gamma$ depends only on $\data$ and will be chosen later.\\
If $\gamma_{**}\ge4\widetilde\gamma$, failure of \eqref{dg-tail} implies
\[
 0\le\ell(t)\le M/4,\qquad M/4\le \mathsf{k}_n(t)\le M .
\]
For each $n$, set
\[
 w_-=(u-\mathsf{k}_n(t))_-,\qquad
 A_n=\{u+\ell(t)<k_n\}\cap Q_n .
\]
On $A_n$, we have $0<w_-\le3M$, and by \eqref{g-lower-level}--\eqref{power-lower-level}, we have
\begin{equation*}
 \begin{gathered}
 c_qM^{q-1}w_-^2\le g_-(u,\mathsf{k}_n(t))
          \le C_qM^{q-1}w_-^2,\\
 (\boldsymbol{k_n(t)}^q-\boldsymbol{ u}^q)_+
          \ge c_qM^{q-1}w_- .
 \end{gathered}
\end{equation*}
Thus, the time-cutoff term in \eqref{energy-full} is bounded by 
\[
 q\iint_{Q_n}g_-(u,k)\partial_t(\zeta^2)\dd x\dd t\leq \frac{C2^nM^{q+1}}{\theta\rho^{2s}}|A_n|.
\]
By \eqref{cutoff-integral}, the spatial cutoff term is estimated by 
\begin{align*}
    \int_{-\theta\rho_n^{2s}}^0\int_{K_n}\int_{K_n}
      &|\zeta(x,t)-\zeta(y,t)|^2
      w_-(x,t)\,w_-(y,t)
      K(x,y,t)\dd y\dd x\dd t\\
&\leq C\,4^n\frac{M^2}{\rho^{2s}}|A_n|.
\end{align*}

Split $\mathbb{R}^N\setminus K_n$ into $K_R\setminus K_n$ and
$\R^N\setminus K_R$. For the first part, the bound $w_-\le3M$ and
\eqref{separation} applied on the support of $\zeta$ give
\[
\int_{-\theta\rho_n^{2s}}^0\int_{K_n}
       \int_{K_R\setminus K_n}
       w_-(y,t)\zeta(x,t)^2w_-(x,t)K(x,y,t)
       \dd y\dd x\dd t
\leq B^n \,C\, \frac{M^{2}}{\rho^{2 s}}\left|A_{n}\right|.
\]
For the second part, using $w_-(y,t)\le M+u_-(y,t)$ and
\eqref{tail-distance} yield,
\begin{align*}
&\int_{-\theta\rho_n^{2s}}^0\int_{K_n}
       \int_{\R^N\setminus K_R}
       w_-(y,t)\zeta(x,t)^2w_-(x,t)K(x,y,t)
       \dd y\dd x\dd t\notag\\
 &\le \frac{C\,M^2}{\rho^{2s}}|A_n|
  +C\int_{-\theta\rho_n^{2s}}^0\int_{K_n}
       \zeta^2w_-
       \int_{\R^N\setminus K_R}
           \frac{u_-(y,t)}{|y|^{N+2s}}\dd y\dd x\dd t .
\end{align*}
Thus, we combine them to get
\begin{align}\label{dg-exterior}
 &\int_{-\theta\rho_n^{2s}}^0\int_{K_n}
       \int_{\R^N\setminus K_n}
       w_-(y,t)\zeta(x,t)^2w_-(x,t)K(x,y,t)
       \dd y\dd x\dd t\notag\\
 &\le \frac{CB^nM^2}{\rho^{2s}}|A_n|
  +C\int_{-\theta\rho_n^{2s}}^0\int_{K_n}
       \zeta^2w_-
       \int_{\R^N\setminus K_R}
           \frac{u_-(y,t)}{|y|^{N+2s}}\dd y\dd x\dd t .
\end{align}
Here and below $B>1$ depends only on $\data$ and may be
enlarged once to dominate all geometric factors.\\
Since $\mathsf{k}_n'(t)=-\ell'(t)$ and $M/4\le \mathsf{k}_n(t)\le M $, by \eqref{power-lower-level} with $a=u(x,t),\,k=\mathsf{k}_n(t)$, we derive
\begin{align*}
\iint_{Q_n} \zeta(x,t)^2 \mathsf{k}_{n}^{\prime}(t)\left(\boldsymbol{u}^{q}-\boldsymbol{\mathsf{ k}_{n}(t)}^{q}\right)_{-} d x d t 
\leq -c_qM^{q-1}\int_{-\theta\rho_n^{2s}}^0
\int_{K_n}\ell'\zeta^2w_-\dd x\dd t.
\end{align*}
Choose $\widetilde\gamma:=\frac{C}{c_q}$ to cancel the last term of
\eqref{dg-exterior}. The same choice works for every $n$.

\smallskip\noindent\textit{Step 2: the energy bound.}
Discarding the nonnegative mixed term in \eqref{energy-full}
and using Remark~\ref{energy-use}, we obtain
\begin{align}\label{dg-energy}
 M^{q-1}&\esssup_{-\theta\rho_n^{2s}<t<0}
        \int_{K_n}w_-^2\zeta^2\dd x+\int_{-\theta\tilde{\rho}_n^{2s}}^0
       \int_{\widetilde{K}_n}\int_{\widetilde{K}_n}
       \frac{|w_-(x,t)-w_-(y,t)|^2}
            {|x-y|^{N+2s}}\dd y\dd x\dd t\notag\\
 &\le \frac{CB^n}{\rho^{2s}}
             \left(\frac{M^{q+1}}{\theta}+M^2\right)|A_n|
  \le\frac{CB^nM^2}{\rho^{2s}}|A_n|.
\end{align}
The last inequality uses $\theta=M^{q-1}$.

\smallskip\noindent\textit{Step 3: embedding and iteration.}
Choose a second cutoff $\varphi$ supported in $\widetilde Q_n$,
equal to one in $\overline Q_n$, with $|D\varphi|\le C2^n/\rho$.
An application of Lemma~\ref{embedding} with $v=w_-\varphi$, $I=(-\theta\tilde{\rho}_n^{2s},0)$, and \eqref{product-cutoff} with $\zeta=\varphi$  gives that
\begin{align}\label{dg-embedding}
 &\iint_{\overline{ Q}_n}w_-^{\widehat m}\dd x\dd t
 \le  \int_{-\theta \tilde{\rho}_{n}^{2s}}^{0} \int_{\widetilde{K}_{n}} | w_{-}\, \varphi(x,t)|^{\widehat{m}} \dd x \dd t\notag\\
 &\qquad\leq C\,\left(\tilde{\rho}_{n}^{2s} \int_{-\theta \tilde{\rho}_{n}^{2s}}^{0} \int_{\widetilde{K}_{n}} \int_{\widetilde{K}_{n}} \frac{\mid w_{-}\, \varphi(x,t)-w_{-}\, \varphi(y, t)|^{2}}{|x-y|^{N+2 s}} \dd y \dd x\right. \notag\\
&\quad\quad\quad\quad\quad\quad\quad\quad\quad\quad\quad\quad \left.+2^{n(N+2s)} \int_{-\theta \tilde{\rho}_{n}^{2s}}^{0} \int_{\widetilde{K}_{n}}\left|w_{-} \varphi(x,t)\right|^{2} \dd x \dd t\right) \notag\\
&\quad\quad\quad\quad\quad \times\left(\esssup_{t\in\left(-\theta \tilde{\rho}_{n}^{2s}, 0\right)}\bint_{\widetilde{K}_n} |w_-\varphi(x, t)|^{2} \dd x \dd t\right)^{\frac{2s}{N}}\notag\\
&\qquad \stackrel{\eqref{product-cutoff}}{\le}  C\, {\left(\rho^{2s} \int_{-\theta \tilde{\rho}_{n}^{2s}}^{0} \int_{\widetilde{K}_{n}} \int_{\widetilde{K}_{n}} \frac{\left|w_{-}\left(x,t\right)-w_{-}(y, t)\right|^{2}}{|x-y|^{N+2s}} \dd y \dd x \dd t\right.} \notag\\
& \quad\quad\quad\left.+B^{n} \int_{-\theta \tilde\rho_{n}^{2 s}}^{0} \int_{\widetilde{K}_{n}} w_{-}^{2}(y, t) \dd y \dd t\right)\times\left(\esssup_{\left.(-\theta\rho_{n}^{2 s}, 0\right)} \bint_{K_{n}}\left|w_{-}(x,t)\right|^{2} \zeta(x, t)^2 \dd x\right)^\frac{2s}{N}\notag\\
 &\qquad \le CB^n M^{\widehat m}
        M^{(1-q)\frac{2s}{N}}\rho^{-2s(1+\frac{2s}{N})}
        |A_n|^{1+\frac{2s}{N}}.
\end{align}
The last inequality is derived by \eqref{dg-energy}.

On $A_{n+1}$, $w_-\ge k_n-k_{n+1}=M/2^{n+2}$.
Using $Q_{n+1}\subset\overline{ Q}_n$, divide
\eqref{dg-embedding} by $(M/2^{n+2})^{\widehat m}$ to get
\[
 |A_{n+1}|
 \le CB^nM^{(1-q)\frac{2s}{N}}\rho^{-2s(1+\frac{2s}{N})}
                    |A_n|^{1+\frac{2s}{N}}.
\]
Set $Y_n=|A_n|/|Q_n|$. The ratios $|Q_n|/|Q_{n+1}|$ are
bounded, and $|Q_n|\asymp\theta\rho^{N+2s}$.
Since
${\frac{2s}{N}}(N+2s)=2s(1+\frac{2s}{N})$ and
$M^{(1-q)\frac{2s}{N}}\theta^{\frac{2s}{N}}=1$, the powers cancel:
\begin{equation*}
 Y_{n+1}\le CB^nY_n^{1+\frac{2s}{N}}.
\end{equation*}
Lemma~\ref{fast} gives $Y_n\to0$ if
$Y_0\le\nu$, with
$\nu=\min\{1/2,C^{-N/(2s)}B^{-N^2/(4s^2)}\}$.
The initial assumption ensures this because
$A_0\subset\{u\le M\}\cap Q_\rho(\theta)$.
It follows that $u+\ell(t)\ge M/2$ almost everywhere in
$Q_{\rho/2}(\theta)$. Since $\ell\le M/4$, this is
\eqref{dg-conclusion}.
\end{proof}

\section{Measure estimates near zero}\label{sec:near}
We recall that the amplitude assumption made at the beginning of Section~\ref{sec:away-dg} is still assumed here.
\subsection{A time slice with positive measure}
\begin{lemma}\label{slice}
Let $\nu\in(0,1)$. If
\[
 |\{u\le M\}\cap Q_\rho (\theta)|
       >\nu |Q_{\rho}(\theta)|,
\]
then there is a time $\tau\in(-\theta\rho^{2s},-\frac{\nu}{2}\theta\rho^{2s})$ such that
\[
 |\{u(\cdot,\tau)\le M\}\cap K_\rho|
       >\frac\nu2|K_\rho|.
\]
The time $\tau$ may be chosen to be an admissible time slice
for the energy estimates. The same statement holds for $-u$.
\end{lemma}
\begin{proof}
If the slice estimate failed almost everywhere on
$(-\theta\rho^{2s},-\frac{\nu}{2}\theta\rho^{2s})$, Fubini's theorem would give
\begin{align*}
 |\{u\le M\}\cap Q_\rho(\theta)|
 &\le\frac\nu2|K_\rho|(1-\nu/2)\theta\rho^{2s}
                       +|K_\rho|\frac\nu2\theta\rho^{2s}\\
 &=(\nu-\nu^2/4)|Q_\rho(\theta)|<\nu|Q_\rho(\theta)|.
\end{align*}
This is a contradiction. The set of qualifying slices has
positive measure, so it intersects the full-measure set of
admissible slices.
\end{proof}

\subsection{Propagation of measure information}
\begin{lemma}\label{timefor}
Let $u$ be a local weak subsolution, let $m\in(0,1]$,
and suppose $2\rho\le R$. If
\[
 |\{u(\cdot,\tau)\le2M-mM\}\cap K_\rho|
       >\frac\nu2|K_\rho|,
\]
then there exist $\delta\in(0,1)$ and
$\varepsilon\in(0,1/2)$, depending only on $\data$ and $\nu$,
such that
\begin{equation}\label{measure-propagation}
 |\{u(\cdot,t)\le2M-\varepsilon mM\}\cap K_\rho|
       \ge\frac\nu4|K_\rho|
\end{equation}
for almost every $t\in(\tau,\tau+\delta\theta\rho^{2s})$,
provided
\begin{equation}\label{propagation-tail}
 \Tail(u_+;Q)\le\delta mM^q .
\end{equation}
The interval must be contained in the reference time interval.
The constants $\delta,\varepsilon$ are independent of $m$.
Both the interval and the tail bound may be shortened by
replacing $\delta$ with a smaller positive number.
For supersolutions, apply the statement to $-u$.
\end{lemma}
\begin{proof}
Put $k=2M-mM$ and $w_+=(u-k)_+$.
Choose a nonnegative spatial cutoff $\zeta$ which equals one
in $K_{(1-\sigma)\rho}$, is supported in
$K_{(1-\sigma/2)\rho}$, and satisfies
$|D\zeta|\le C/(\sigma\rho)$.
The parameter $\sigma\in(0,1/2)$ will be fixed below.
Since $k\ge M$ and $u\le2M$, the active values of $u$ lie in $[M,2M]$.

For a terminal time $t\le\tau+\delta\theta\rho^{2s}$,
Proposition~\ref{energy} gives
\begin{equation}
 \begin{aligned}
 &q\int_{K_\rho}g_+(u(t),k)\zeta^2\dd x-q\int_{K_\rho}g_+(u(\tau),k)\zeta^2\dd x\\
  &\le \int_{\tau}^{\tau+\delta\theta\rho^{2s}}\int_{K_\rho}\int_{K_\rho}
      |\zeta(x)-\zeta(y)|^2
      w_+(x,t)\,w_+(y,t)
      K(x,y,t)\dd y\dd x\dd t
 \notag\\
&\quad+\int_{\tau}^{\tau+\delta\theta\rho^{2s}}\int_{K_\rho}\int_{\R^N\setminus K_\rho}
      w_+(y,t)\zeta(x)^2w_+(x,t)
      K(x,y,t)\dd y\dd x\dd t .
 \end{aligned}
\end{equation}
Since $w_+(x,t)\leq mM$ for $x\in K_R$, and  \eqref{cutoff-integral}, the first term on the right-hand side is estimated by
\begin{align*}
&\int_{\tau}^{\tau+\delta\theta\rho^{2s}}\!\!\int_{K_{\rho}}\!\int_{K_{\rho}}
\bigl(\zeta(x)-\zeta(y)\bigr)^2 w_+(y,t)w_+(x,t)K(x,y,t)\,\dd x\dd y\dd t\notag\\
&\quad\leq C\,\frac{(mM)^2\theta\delta}{\sigma^2}\,\meas{K_{\rho}}.
\end{align*}
Regarding the second term, for $y\notin K_R$, $w_+(y,t)\le u_+(y,t)$ because $k>0$; for $y\in K_R\setminus K_\rho,\, x\in K_{(1-\sigma/2)\rho}$, we use $w_+(\cdot,t)\leq mM$, thus

\begin{align}\label{propogation-interaction}
&\int_{\tau}^{\tau+\delta\theta\rho^{2s}}\!\!\int_{\mathbb R^N\setminus K_{\rho}}\!\int_{K_{\rho}}
w_+(y,t)\,\zeta(x)^2\,w_+(x,t)\,K(x,y,t)\,\dd x\dd y\dd t\notag\\
&\quad\le mM\int_{\tau}^{\tau+\delta\theta\rho^{2s}}\!\!\int_{K_{\rho}}\zeta(x)^2\int_{\mathbb R^N\setminus K_R}\frac{u_+(y,t)}{|x-y|^{N+2s}}\,\dd y\dd x\dd t\notag\\
&\quad\quad+\,(mM)^2\int_{\tau}^{\tau+\delta\theta\rho^{2s}}\!\!\int_{K_{\rho}}\zeta(x)^2\int_{K_R\setminus K_{\rho}}\frac{1}{|x-y|^{N+2s}}\,\dd y \dd x\dd t\notag\\
& \stackrel{\eqref{tail-distance},\eqref{separation} }{\leq} C\, mM\,\meas{K_{\rho}}\Tail(u_+;Q)
+C\,\frac{(mM)^2\,\meas{K_{\rho}}\delta\theta}{\sigma^{N+2s}}
\end{align}

Thus, combining above estimates with $\theta=M^{q-1}$ and \eqref{propagation-tail}, we have for all $\tau<t\leq \tau+\delta\theta\rho^{2s}$,that
\begin{equation}\label{prop-energy}
 \begin{aligned}
 &q\int_{K_\rho}g_+(u(t),k)\zeta^2\dd x-q\int_{K_\rho}g_+(u(\tau),k)\zeta^2\dd x\\
 &\qquad\le C\delta M^{q-1}(mM)^2
        (1+\sigma^{-2}+\sigma^{-N-2s})|K_\rho|.
 \end{aligned}
\end{equation}

By initial measure assumption, the second term on left-hand side of \eqref{prop-energy} is bounded by
\begin{align}\label{exact-time-energy01}
q\int_{K_\rho}g_+(u(\tau),k)\zeta^2\dd x\leq qg_+(2M,k)(1-\nu/2)|K_\rho|.
\end{align}

For the first term on the left-hand side of \eqref{prop-energy}, use $\zeta=1$ in $K_{(1-\sigma)\rho}$
and the monotonicity of $a\mapsto g_+(a,k)$ for $a\ge k$:
\begin{align*}
 q\int_{K_\rho}g_+(u(t),k)\zeta^2\dd x
 &\ge q\int_{K_{(1-\sigma)\rho}\cap
                 \{u(t)>2M-\varepsilon mM\}}
                 g_+(u(t),k)\dd x\notag\\
 &\ge qg_+(2M-\varepsilon mM,k)
       |\{u(t)>2M-\varepsilon mM\}\cap K_{(1-\sigma)\rho}|.
\end{align*}
Since $k=2M-mM\ge M$ and $\varepsilon\le1/2$, we estimate
\begin{align*}
 g_+(2M-\varepsilon mM,k)
 &=\int_k^{2M-\varepsilon mM}r^{q-1}(r-k)\dd r\notag\\
 &\ge c_qM^{q-1}(mM)^2, \\
 g_+(2M,k)-g_+(2M-\varepsilon mM,k)
 &=\int_{2M-\varepsilon mM}^{2M}r^{q-1}(r-k)\dd r\notag\\
 &\le C_qM^{q-1}(mM)(\varepsilon mM). 
\end{align*}
Consequently, we get
\begin{align}\label{exact-time-energy05}
 \frac{g_+(2M,k)}{g_+(2M-\varepsilon mM,k)}
       \le1+C_q\varepsilon
\end{align}
uniformly in $m$.
Combining these estimates with
\eqref{prop-energy} and
$|K_\rho\setminus K_{(1-\sigma)\rho}|\le N\sigma|K_\rho|$
yields
\begin{align}\label{prop-ratio}
 &\frac{|\{u(t)>2M-\varepsilon mM\}\cap K_\rho|}{|K_\rho|}\\\notag
 &\qquad\le (1+C_q\varepsilon)(1-\nu/2) +C\delta(1+\sigma^{-2}+\sigma^{-N-2s})+N\sigma .
\end{align}
Choose first $\sigma=\nu/(16N)$, then $\varepsilon$ so that
$C_q\varepsilon(1-\nu/2)\le\nu/8$, and finally $\delta$
so that
$C\delta(1+\sigma^{-2}+\sigma^{-N-2s})\le\nu/16$.
The right-hand side of \eqref{prop-ratio} is at most
$1-\nu/4$. Taking the complement proves
\eqref{measure-propagation}.
The same estimates hold on every shorter time interval.
\end{proof}

\subsection{A shrinking lemma}
\begin{lemma}\label{shrink}
Let $u$ be a locally bounded local weak subsolution in $Q$.
Let $\theta=M^{q-1}$, $\sigma,\varepsilon,\delta\in(0,1/2)$,
and $2\rho\le R$. Write
$I=(\tau,\tau+\delta\theta\rho^{2s})$ and assume
$K_{2\rho}\times I\subset Q$. If
\[
 |\{u(\cdot,t)\le2M-\varepsilon M\}\cap K_\rho|
       \ge\frac\nu4|K_\rho|\quad\text{for a.e. }t\in I,
\]
then either
\begin{equation}\label{shrinking-tail}
 \Tail(u_+;Q)>\frac{\sigma\varepsilon}{2}M^q
\end{equation}
or
\begin{equation}\label{shrinking-conclusion}
 \frac{|\{u\ge2M-\sigma\varepsilon M/2\}\cap(K_\rho\times I)|}
      {|K_\rho\times I|}
 \le\frac{C\sigma}{\nu}(1+\delta^{-1}),
\end{equation}
where $C$ depends only on $\data$.
For supersolutions, replace $u$ with $-u$.
\end{lemma}
\begin{proof}
Assume that \eqref{shrinking-tail} fails.
Take $k=2M-\sigma\varepsilon M$, set $w_\pm=(u-k)_\pm$,
and choose a spatial cutoff
$\zeta=1$ in $K_\rho$, supported in $K_{3\rho/2}$, with
$|D\zeta|\le C/\rho$.
Apply \eqref{energy-full} in $K_{2\rho}\times I$ and keep the
 term involving mixed truncations: 
\begin{align*}
&\int_I\int_{K_{2\rho}}\int_{K_{2\rho}}
       w_-(y,t)w_+(x,t)\zeta(x)^2K(x,y,t)\dd y\dd x\dd t
\leq q\int_{K_{2\rho}}g_+(u(\tau),k)\,\zeta(x)^2\,\dd x\\
&\qquad+\int_{I}\int_{K_{2\rho}}\int_{K_{2\rho}}
      |\zeta(x)-\zeta(y)|^2
      w_+(x,t)\,w_+(y,t)
      K(x,y,t)\dd y\dd x\dd t
 \notag\\
 &\qquad+\int_{I}\int_{K_{2\rho}}\int_{\R^N\setminus K_{2\rho}}
      w_+(y,t)\zeta(x)^2w_+(x,t)
      K(x,y,t)\dd y\dd x\dd t.
\end{align*}

Observe that
$w_+\le\sigma\varepsilon M$ in $Q$. By \eqref{g-lower-level}, the first term on the right-hand side is estimated by
\[
q\int_{K_{2\rho}}g_+(u(\tau),k)\,\zeta(x)^2\,\dd x
\leq CM^{q-1}(\sigma\varepsilon M)^2|K_\rho|.
\]
As for the second term, \eqref{cutoff-integral} gives
\[
\int_{I}\int_{K_{2\rho}}\int_{K_{2\rho}}
      |\zeta(x)-\zeta(y)|^2
      w_+(x,t)\,w_+(y,t)
      K(x,y,t)\dd y\dd x\dd t
\leq C(\sigma\varepsilon M)^2\rho^{-2s}|K_\rho\times I|.
\]
The third term can be estimated by computations similar to those in \eqref{propogation-interaction}
\[
\begin{aligned}
&\int_{I}\int_{K_{2\rho}}\int_{\R^N\setminus K_{2\rho}}
      w_+(y,t)\zeta(x)^2w_+(x,t)
      K(x,y,t)\dd y\dd x\dd t\\
 &\qquad\leq C\, \sigma\varepsilon M\,\meas{K_{\rho}}\Tail(u_+;Q)
+C(\sigma\varepsilon M)^2\rho^{-2s}|K_\rho\times I|.
\end{aligned}
\]
Recalling $\Tail(u_+;Q)\leq\frac{\sigma\varepsilon}{2}M^q$, then combining all these estimates, we get
\begin{align}\label{shrink-mixed-upper}
 &\int_I\int_{K_{2\rho}}\int_{K_{2\rho}}
       w_-(y,t)w_+(x,t)\zeta(x)^2K(x,y,t)\dd y\dd x\dd t
 \notag\\
 &\qquad\le
 \frac{C(\sigma\varepsilon M)^2}{\rho^{2s}}
              (1+\delta^{-1})|K_\rho\times I|.
\end{align}

For the left-hand side of \eqref{shrink-mixed-upper}, if $u(y,t)\le2M-\varepsilon M$ and
$u(x,t)\ge2M-\sigma\varepsilon M/2$, then
\[
 w_-(y,t)\ge(1-\sigma)\varepsilon M,\qquad
 w_+(x,t)\ge\sigma\varepsilon M/2.
\]
For $x,y\in K_\rho$, we also have $|x-y|\le2\rho$.
Using the lower bound on the measure of
$\{u(\cdot,t)\le2M-\varepsilon M\}\cap K_\rho$, we obtain
\begin{align}\label{shrink-mixed-lower}
 &\int_I\int_{K_{2\rho}}\int_{K_{2\rho}}
       w_-(y,t)w_+(x,t)
       \zeta(x)^2 K(x,y,t)\dd y\dd x\dd t
       \notag\\
 &\quad\ge
 \frac{C\,\nu\sigma(1-\sigma)(\varepsilon M)^2}{\rho^{2s}}
 \int_I|\{u(\cdot,t)\ge2M-\sigma\varepsilon M/2\}
                \cap K_\rho|\dd t .
\end{align}
Combine \eqref{shrink-mixed-upper} and
\eqref{shrink-mixed-lower}, and use $1-\sigma\ge1/2$.
This proves \eqref{shrinking-conclusion}.
\end{proof}

\subsection{A De Giorgi lemma near the upper bound}
\begin{lemma}\label{De03}
Let $u$ be a locally bounded local weak subsolution in $Q$.
Let $\delta,\xi\in(0,1)$ and $\theta=M^{q-1}$.
Assume $2\rho\le R$ and
$K_\rho\times(\tau,\tau+\delta\theta\rho^{2s})\subset Q$.
After enlarging $\gamma_{**}$ in Lemma~\ref{De01} if necessary,
there is
$\nu_0=\nu_0(\delta,\data)\in(0,1)$, independent of $\xi$,
such that
\begin{align}
 |\{u\ge2M-\xi M\}\cap(K_\rho\times(\tau,\tau+\delta\theta\rho^{2s}))|
 &\le\nu_0|K_\rho|\delta\theta\rho^{2s},\label{upper-smallness}\\
 \gamma_{**}\Tail(u_+;Q)&\le\xi M^q\label{upper-tail}
\end{align}
imply
\begin{equation*}
 u\le2M-\frac{\xi M}{4}
 \quad\text{a.e. in }
K_{\rho/2}\times(\tau+3\delta\theta\rho^{2s}/4,\tau+\delta\theta\rho^{2s}).
\end{equation*}
The analogous assertion holds for $-u$.
\end{lemma}
\begin{proof}
Translate the final time to zero and keep the radii
$\rho_n,\hat\rho_n,\tilde\rho_n,\bar\rho_n$ and the balls $K_{n}, \widehat{K}_{n},\widetilde{K}_{n}, \overline{K}_{n} $ used in
Lemma~\ref{De01}. To obtain exactly the final quarter of the
time interval, define the heights independently of the radii:
\[
 \begin{gathered}
 d_n=\delta\theta\rho^{2s}\left(\frac14+\frac{3}{4\cdot2^n}\right),\\
 \hat d_n=\frac{3d_n+d_{n+1}}4,\qquad
 \tilde d_n=\frac{d_n+d_{n+1}}2,\qquad
 \bar d_n=\frac{d_n+3d_{n+1}}4 .
 \end{gathered}
\]
Set $$
\left\{\begin{array}{l}
Q_{n}=K_n \times\left(-d_n, 0\right), \quad \widetilde{Q}_{n}=\widetilde{K}_{n} \times\left(-\tilde{d}_n, 0\right), \\
\widehat{Q}_{n}=\widehat{K}_{n} \times\left(-\hat{d}_n, 0\right), \quad \overline{Q}_{n}=\overline{K}_{n} \times\left(-\bar{d}_n, 0\right).
\end{array}\right.
$$

Choose $\zeta$ supported in $\widehat Q_n$, equal to one in
$\widetilde Q_n$, and satisfying
\[
 |D\zeta|\le C2^n/\rho,\qquad
 |\partial_t\zeta|\le C2^n/\delta\theta\rho^{2s} .
\]
Support here concerns the lateral boundary and the initial
time; the cutoff need not vanish at the final time.

Take
\begin{equation*}
 k_n=2M-\frac{\xi M}{2}-\frac{\xi M}{2^{n+1}},
 \qquad \mathsf{k}_n(t)=k_n+\ell(t),
\end{equation*}
where
\[
 \ell(t)=\widetilde\gamma M^{1-q}\int_{T_1}^t
       \int_{\R^N\setminus K_R}
       \frac{u_+(y,\tau)}{|y|^{N+2s}}\dd y\dd\tau.
\]
Choose $\gamma_{**}\ge4\widetilde\gamma$. Assumption \eqref{upper-tail}
then gives $0\le\ell\le\xi M/4$.
In particular $\mathsf{k}_n(t)\in [ M,2M]$.
Let $w_+=(u-\mathsf{k}_n(t))_+$ and $A_n=\{u>\mathsf{k}_n(t)\}\cap Q_n$.
Then $w_+\le\xi M$ in $Q_n$ and
\[
 \begin{gathered}
 c_qM^{q-1}w_+^2\le g_+(u,\mathsf{k}_n(t))
                      \le C_qM^{q-1}w_+^2,\\
 (\boldsymbol{u}^q-\boldsymbol{ \mathsf{k}_n(t)}^q)_+
                      \ge c_qM^{q-1}w_+.
 \end{gathered}
\]

We spell out the energy estimates to keep track of the
dependence on $\delta$ and $\xi$.
As in Lemma~\ref{De01}, $B>1$ is a fixed constant depending
only on the $\data$, chosen to dominate the powers of $2$
from the cutoff bounds and the level differences. In
particular, $B$ is independent of $n$, $\delta$, and $\xi$.
Then ,the time-cutoff term on the left-hand side of \eqref{energy-full} is bounded by
\[
q\iint_{Q_n}g_+(u,k)\partial_t(\zeta^2)\dd x\dd t\leq
CB^n M^{q-1}(\xi M)^2(\delta\theta\rho^{2s})^{-1}|A_n|.
\]

The spatial cutoff term is estimated by 
\[
\int_{-d_n}^{0}\int_{K_n}\int_{K_n}
      |\zeta(x,t)-\zeta(y,t)|^2
      w_+(x,t)\,w_+(y,t)
      K(x,y,t)\dd y\dd x\dd t
      \leq CB^n\frac{(\xi M)^2}{\rho^{2s}}|A_n|.
\]
For the exterior interaction, $w_+(y,t)\le\xi M$ in
$K_R\setminus K_n$ and $w_+(y,t)\le u_+(y,t)$ outside $K_R$.
Thus it is bounded by
\begin{align*}
&\int_{-d_n}^0\int_{K_n}
       \int_{\R^N\setminus K_n}
       w_+(y,t)\zeta(x,t)^2w_+(x,t)K(x,y,t)
       \dd y\dd x\dd t\\
 &\leq \frac{CB^n(\xi M)^2}{\rho^{2s}}|A_n|
 +C\int_{-d_n}^0\int_{K_n}\zeta^2w_+(x,t)
       \int_{\R^N\setminus K_R}
        \frac{u_+(y,t)}{|y|^{N+2s}}\dd y\dd x\dd t .
\end{align*}
The last term is absorbed by the term
$-c_qM^{q-1}\int\ell'\zeta^2w_+$, derived by the second term on the left-hand side of \eqref{energy-full},
after choosing $\widetilde\gamma=\frac{C}{c_q}$.
Consequently,
\begin{align}\label{upper-energy}
M^{q-1} &\esssup_{-d_n<t<0}\int_{K_n}w_+^2\zeta^2\dd x
+\int_{-\tilde{d}_n}^0\int_{\widetilde{K}_n}\int_{\widetilde{K}_n}
   \frac{|w_+(x,t)-w_+(y,t)|^2}{|x-y|^{N+2s}}
                   \dd y\dd x\dd t\notag\\
 &\le\frac{CB^n(\xi M)^2}{\rho^{2s}}
                      (1+\delta^{-1})|A_n|.
\end{align}

Apply Lemma~\ref{embedding} with the second cutoff supported
in $\widetilde Q_n$ and equal to one in  $\overline{Q}_n$, just as in
Lemma~\ref{De01}. By arguments similar to those in Step 3 of Lemma~\ref{De01}, we get
\begin{align*}
 &\iint_{\overline{Q}_n}w_+^{\widehat m}\dd x\dd t\\
&\quad{\le}  C\, {\left(\rho^{2s} \int_{-\tilde{d}_n}^{0} \int_{\widetilde{K}_{n}} \int_{\widetilde{K}_{n}} \frac{\left|w_{+}\left(x,t\right)-w_{+}(y, t)\right|^{2}}{|x-y|^{N+2s}} \dd y \dd x \dd t\right.} \notag\\
& \quad\quad\quad\left.+B^{n} \int_{-\tilde{d}_n}^{0} \int_{\widetilde{K}_{n}} w_{+}^{2}(y, t) \dd y \dd t\right)\times\left(\esssup_{\left.(-d_n, 0\right)} \bint_{K_{n}}\left|w_{+}(x,t)\right|^{2} \zeta(x, t)^2 \dd x\right)^\frac{2s}{N}\notag\\
 &\qquad \le CB^n {(\xi M)}^{\widehat m}
        M^{(1-q)\frac{2s}{N}}\rho^{-2s(1+\frac{2s}{N})}(1+\delta^{-1})^{1+\frac{2s}{N}}
        |A_n|^{1+\frac{2s}{N}}.
\end{align*}
The last inequality is derived by \eqref{upper-energy}.

Since 
$w_+\geq k_{n+1}-k_n=\xi M/2^{n+2}$ on $A_{n+1}$, one obtains
\[
 |A_{n+1}|\le CB^n M^{(1-q)\frac{2s}{N}}
       \rho^{-2s(1+\frac{2s}{N})}(1+\delta^{-1})^{1+\frac{2s}{N}}
       |A_n|^{1+\frac{2s}{N}}.
\]

Normalizing by $|Q_n|\asymp
\delta M^{q-1}\rho^{N+2s}$ yields
\begin{equation*}
 Y_{n+1}\le
 C B^n\delta^{\frac{2s}{N}}(1+\delta^{-1})^{1+\frac{2s}{N}}
           Y_n^{1+\frac{2s}{N}},\qquad Y_n=\frac{|A_n|}{|Q_n|}.
\end{equation*}
It suffices to set
\begin{equation}\label{nu0}
 \nu_0(\delta)=
 \min\left\{\frac12,\,
 \Big[C\delta^{\frac{2s}{N}}(1+\delta^{-1})^{1+\frac{2s}{N}}\Big]^{-\frac{N}{2s}}
                B^{-\frac{N^2}{4s^2}}\right\}.
\end{equation}
Since  $k_0=2M-\xi M$ and $\ell\ge0$,
\eqref{upper-smallness} ensures $Y_0\le\nu_0(\delta)$.
The iteration gives
$u\le2M-\xi M/2+\ell(t)\le2M-\xi M/4$
in the limiting cylinder, as required.
\end{proof}

\subsection{Reaching the final time by finitely many steps}
\begin{lemma}\label{EoP}
Let $u$ be a locally bounded local weak subsolution in $Q$.
Let $\nu$ be the constant in Lemma~\ref{De01}, and assume
$K_{2\rho}\times(-\theta\rho^{2s},0)\subset Q$.
There exist $\eta\in(0,1/2)$ and $\gamma\in(0,1/4)$,
depending only on $\data$, such that
\[
 |\{u\le M\}\cap(K_\rho\times(-\theta\rho^{2s},0))|
     >\nu|K_\rho|\theta\rho^{2s},\qquad
 \Tail(u_+;Q)\le\gamma M^q
\]
imply
\begin{equation}\label{EoP-conclusion}
 u\le2M(1-\eta)\quad\text{a.e. in }
 K_{\rho/2}\times(-\tfrac{7}{16}\nu\theta\rho^{2s},0).
\end{equation}
For a supersolution, apply the statement to $-u$.
\end{lemma}
\begin{proof}
We combine the time-slice, propagation, shrinking, and
De Giorgi estimates in a finite sequence.
The reference amplitude $M$ and the scale $\theta=M^{q-1}$
remain fixed throughout the proof.

\smallskip\noindent\textit{Step 1: the first interval.}
Lemma~\ref{slice} gives an admissible
$\tau\in(-\theta\rho^{2s},-\nu \theta\rho^{2s}/2)$ with
$|\{u(\tau)\le M\}\cap K_\rho|>\nu|K_\rho|/2$.
Use Lemma~\ref{timefor} with $m=1$ and 
$\delta, \varepsilon$, depending on $\nu$ and data. Note that $\delta\le\nu/16$.
Let $\varepsilon$ be its propagation parameter.

By Lemma~\ref{shrink}, we have
\[
\bigl|\{u\ge 2M-\frac{\sigma\varepsilon}{2} M\}\cap K_{\rho}\times(\tau,\tau+\delta\theta\rho^{2s})\bigr|\le\frac{C\sigma(1+\frac1\delta)}{\nu}\bigl|K_{\rho}\times(\tau,\tau+\delta\theta\rho^{2s})\bigr|,
\]
where $C$ is the
structural constant in that lemma.
Choose $\sigma\in(0,1/2)$ so small that
\[
 C\sigma(1+\delta^{-1})\nu^{-1}\le\nu_0(\delta).
\]
This is a noncircular choice because \eqref{nu0} does not
depend on $\xi=\sigma\varepsilon/2$.
Lemma~\ref{De03} gives
\begin{equation}\label{EoP-first}
 u\le2M-m_1M
 \quad\text{in }K_{\rho/2}\times(a_0,b_0),\qquad
 m_1=\frac{\sigma\varepsilon}{8}.
\end{equation}
Thus $m_1M$ is the first gap below the upper bound $2M$.
Here
\[
 a_0=\tau+\frac34\delta \theta\rho^{2s},\qquad b_0=\tau+\delta \theta\rho^{2s}.
\]
The required tail conditions are all fulfilled if
\begin{equation*}
 \Tail(u_+;Q)\le
 \min\{\delta,\sigma\varepsilon/(2\gamma_{**})\}M^q.
\end{equation*}
Here $\gamma_{**}\ge1$ also covers the shrinking condition.
Moreover,
\[
 a_0\le-\frac\nu2 \theta\rho^{2s}+\frac{3\nu}{64}\theta\rho^{2s}
       =-\frac{29\nu}{64}\theta\rho^{2s}<-\frac{7\nu}{16}\theta\rho^{2s} .
\]
If $b_0\geq 0$, the proof is finished; otherwise, we proceed the proof.\\
\smallskip\noindent\textit{Step 2: a step with uniform parameters.}
Fix $\bar\nu=2^{-N}$.
Let $m\in(0,1]$ and suppose that, at an admissible time $s$,
$u(\cdot,s)\le2M-mM$ in $K_{\rho/2}$.
The number $mM$ is the current gap below $2M$; on the
first interval we have $m=m_1$ from \eqref{EoP-first}.
This bound gives
\[
\bigl|\{u(x,s)\le 2M-m M\}\cap K_{\rho}\bigr|>\bigl|\{u(x,s)\le 2M-m M\}\cap K_{\frac{1}{2}\rho}\bigr|>\frac12|K_{\frac{\rho}{2}}|=\frac{\bar{\nu}}{2}|K_{\rho}|.
\]

Apply Lemma~\ref{timefor} with $\bar\nu$ in place of $\nu$.
Write its constants as $\bar\varepsilon,\bar\delta$ and
decrease $\bar\delta$ so that $\bar\delta\le\delta/4$.
These constants are independent of $m$.

For later use we permit any duration $d\theta \rho^{2s}$ with
$d\in[\bar\delta/2,\bar\delta]$.
The expression in \eqref{nu0} has a positive minimum on this
interval. Choose $\bar\sigma\in(0,1/2)$ so that
\begin{equation}\label{EoP-uniform-choice}
 \frac{C\bar\sigma(1+2/\bar\delta)}{\bar\nu}
 \le\inf_{\bar\delta/2\le d\le\bar\delta}\nu_0(d).
\end{equation}
Starting with depth $m$, the three lemmas now give depth
$am$ on the final quarter of the new interval, where
\[
 a=\frac{\bar\sigma\bar\varepsilon}{8}\in(0,1).
\]
It is sufficient that
\begin{equation*}
 \Tail(u_+;Q)
 \le c_*mM^q,\qquad
 c_*=\min\{\bar\delta/2,\bar\sigma\bar\varepsilon/(2\gamma_{**})\}.
\end{equation*}
Indeed, the propagated depth is $\bar\varepsilon m$, and
the level parameter for Lemma~\ref{De03} is
$\xi=\bar\sigma\bar\varepsilon m/2$.

\smallskip\noindent\textit{Step 3: covering the remaining time interval.}
Set $m^{(0)}=m_1=\sigma\varepsilon/8$, as in \eqref{EoP-first}. Suppose that
$u\le2M-m^{(j)}M$ on $K_{\rho/2}\times(a_0,b_j)$, where $b_j<0$.
Put $h=\bar\delta \theta\rho^{2s}$ and choose an admissible starting time
\[
 s_j\in(b_j-3h/4,\ b_j-11h/16).
\]
This interval lies in $(a_0,b_j)$: initially
$b_0-a_0=\delta \theta\rho^{2s}/4\ge h$, and the covered interval only grows.
First consider duration
$d_j\theta\rho^{2s}=\frac43(b_j-s_j)$, which lies in $[11h/12,h]$.
Then its final quarter starts exactly at $b_j$ and ends at
\[
 b_{j+1}=b_j+\frac13(b_j-s_j)
       \ge b_j+\frac{11h}{48}.
\]
If this end would exceed zero, the proof is finished.
All these durations belong to $[\bar\delta/2,\bar\delta]\theta \rho^{2s}$,
as required in \eqref{EoP-uniform-choice}.

At each step set $m^{(j+1)}=am^{(j)}$.
The previous estimate also holds with this smaller depth, so
the union remains bounded by $2M-m^{(j+1)}M$.
Because $b_0>-\theta\rho^{2s}$, the construction reaches zero after at most
\[
 J=\left\lceil\frac{48}{11\bar\delta}\right\rceil+1
\]
additional steps. Define $m_{\mathrm{fin}}=m_1a^J>0$.
All the tail conditions are simultaneously satisfied if
\[
 \gamma\le
 \min\{1/8,\delta,\sigma\varepsilon/(2\gamma_{**}),c_*m_{\mathrm{fin}}\}.
\]
Finally take $\eta=m_{\mathrm{fin}}/2$. The estimate holds on
$K_{\rho/2}\times(a_0,0)$, which contains the cylinder in
\eqref{EoP-conclusion}.
\end{proof}

\section{Proof of Theorem~\ref{main}}\label{sec:continuity}
We first carry out the intrinsic iteration near zero. If the
solution separates from zero at a finite step, we normalize
it and apply the oscillation argument of \cite{Liao-24b}.

Translate $(x_o,t_o)$ to $(0,0)$ and fix
$Q=K_R\times(T_1,T_2)\Subset E_T$ containing a neighborhood
of that point. Choose $M_0>0$ so that
\begin{equation}\label{initial-amplitude}
 2M_0\ge\esssup_Q|u|
       +2(4\gamma_{**})^{1/q}\Tail(u;Q)^{1/q},
\end{equation}
where $\gamma_{**}$ is the constant in Lemma~\ref{De01}.
This gives both $\esssup_Q|u|\le2M_0$ and
$\Tail(u;Q)\le M_0^q/(4\gamma_{**})$.
The factor in front of the tail is absorbed into the
structural constant of \eqref{main-amplitude}.
Set $\theta_0=M_0^{q-1}$ and choose $\rho_0>0$ such that
$Q_0:=Q_{2\rho_0}(\theta_0)\subset Q$.
For $j\ge0$, we shall use the notation
\[
 \rho_j=\lambda^j\rho_0,\qquad
 \theta_j=M_j^{q-1},\qquad Q_j=Q_{2\rho_j}(\theta_j),\qquad
 \widehat Q_j=K_{2\rho_j}\times(-\theta_j(c\rho_j)^{2s},0).
\]
The constants $c,\lambda\in(0,1)$ are chosen below.
Let the $\nu$ be the constant in Lemma \ref{De01},
the iteration continues as long as
\begin{equation}\label{near-measure}
 \begin{split}
 |\{u\le M_j\}\cap Q_{c\rho_j}(\theta_j)|
       &>\nu|Q_{c\rho_j}(\theta_j)| \text{ and }\\
 |\{-u\le M_j\}\cap Q_{c\rho_j}(\theta_j)|
       &>\nu|Q_{c\rho_j}(\theta_j)|.
 \end{split}
\end{equation}
The first index at which one of these conditions fails is
treated in Section~\ref{sec:away-normalization}.

\subsection{Reduction of oscillation near zero}
\subsubsection{The first step}
Assume \eqref{near-measure} for $j=0$.
Apply Lemma~\ref{EoP} to $u$ and $-u$, with local radius
$c\rho_0$ and reference cylinder $\widehat Q_0$.
The required inclusion is 
\[
 K_{2c\rho_0}\times(-\theta_0(c\rho_0)^{2s},0)
       \subset\widehat Q_0.
\]
If $\Tail(u;\widehat Q_0)\le\gamma M_0^q$, the lemma gives
$|u|\le2M_0(1-\eta)$ in
$K_{c\rho_0/2}\times(-\frac{7}{16}\nu\theta_0(c\rho_0)^{2s},0)$.
By \eqref{initial-amplitude}, and taking
\begin{equation}\label{bar-gamma}
 \overline\gamma\ge
 \max\{2\gamma^{-1/q},\,2(4\gamma_{**})^{1/q},\,2\},
\end{equation}
we obtain
\begin{equation*}
 \begin{gathered}
 \esssup_{K_{c\rho_0/2}\times
       (-\frac{7}{16}\nu\theta_0(c\rho_0)^{2s},0)}|u|
       \le2M_1,\\
 2M_1=\max\{2M_0(1-\eta),\,
                   \overline\gamma\Tail(u;\widehat Q_0)^{1/q}\}.
 \end{gathered}
\end{equation*}
Splitting the tail at radius $R$ and using $|u|\le2M_0$ in $Q$,
\begin{equation}\label{first-tail}
\begin{aligned}
 \Tail(u;\widehat Q_0)&\le
 \int_{T_1}^{T_2}\!\int_{\R^N\setminus K_R}\frac{|u(y,t)|}{|y|^{N+2s}}\,\dd y\dd t+(2M_0)\int_{-(c\rho_0)^{2s}\theta_0}^{0}\!\int_{K_R\setminus K_{2\rho_0}}\frac{1}{|y|^{N+2s}}\,\dd y\dd t\\
 &\leq \Tail(u;Q)+Cc^{2s}M_0^q.
 \end{aligned}
\end{equation}
Consequently, there is $\gamma_*\ge1$, depending only on
the $\data$, such that
\[
 (1-\eta)M_0\le M_1\le\gamma_*M_0.
\]
We fix $$\gamma_*=\max\left\{\overline\gamma \left(\frac{1}{(4\gamma_{**})^\frac{1}{q}}+C\right), 2+ \frac{C}{(1-\eta)^\frac{1}{q}}  \right\}$$
and put
\[
 \widehat\gamma=\max\{\gamma_*^{q-1},(1-\eta)^{q-1}\}.
\]
Then $\theta_1\le\widehat\gamma\theta_0$.
Choose
\begin{equation}\label{lambdaj}
 \lambda\le\min\left\{\frac c4,\,
       \frac c2\left(\frac{7\nu}{32\widehat\gamma}\right)^{1/(2s)}
                 \right\}.
\end{equation}
This ensures
\[
 Q_1\subset K_{c\rho_0/2}\times
       (-\tfrac{7}{16}\nu\theta_0(c\rho_0)^{2s},0)\subset Q_0,
\]
and hence
\begin{equation*}
 \esssup_{Q_1}|u|\le2M_1.
\end{equation*}

\subsubsection{The second step and the induction}
If \eqref{near-measure} holds at $j=1$, repeat the argument in
$\widehat Q_1$ and set
\[
 2M_2=\max\{2M_1(1-\eta),\,
                  \overline\gamma\Tail(u;\widehat Q_1)^{1/q}\}.
\]
The nesting of the time intervals and $|u|\le2M_0$ in $Q_0$ give
\begin{equation*}
\begin{aligned}
 \Tail(u;\widehat Q_1)
 &=\int_{-\theta_1(c\rho_1)^{2s}}^{0}\!\int_{\R^N\setminus K_{2\rho_0}}+\int_{-\theta_1(c\rho_1)^{2s}}^{0}\!\int_{K_{2\rho_0}\setminus K_{2\rho_1}}\frac{|u(y,t)|}{|y|^{N+2s}}\,\dd y \dd t\\
 &\le\Tail(u;\widehat Q_0)+Cc^{2s}M_1^{q-1}M_0 .
 \end{aligned}
\end{equation*}
Since
$\overline\gamma^q\Tail(u;\widehat Q_0)\le(2M_1)^q$ and
$M_0\le M_1/(1-\eta)$, thus
\[
(1-\eta) M_1\leq M_2\leq \gamma_* M_1.
\]

At a general step, define
\begin{equation}\label{recursive}
 2M_{j+1}=
 \max\{2M_j(1-\eta),\,\overline\gamma\Tail(u;\widehat Q_j)^{1/q}\}.
\end{equation}
Lemma~\ref{EoP} bounds $|u|$ by $2M_{j+1}$ on the corresponding
smaller cylinder. For $j\ge1$, the same splitting gives
\begin{align}\label{tailj}
 \Tail(u;\widehat Q_j)
 &\le\Tail(u;\widehat Q_{j-1})
                 +Cc^{2s}M_j^{q-1}M_{j-1}\notag\\
 &\le\left(\frac{2M_j}{\overline\gamma}\right)^q
                 +\frac{Cc^{2s}}{1-\eta}M_j^q.
\end{align}
Inductively, we  have
\begin{equation}\label{induction-bounds}
 \begin{gathered}
 (1-\eta)M_j\le M_{j+1}\le\gamma_*M_j,\qquad
 \theta_{j+1}\le\widehat\gamma\theta_j,\\
 Q_{j+1}\subset K_{c\rho_j/2}\times
       (-\tfrac{7}{16}\nu\theta_j(c\rho_j)^{2s},0)\subset Q_j,
 \qquad \esssup_{Q_{j+1}}|u|\le2M_{j+1}.
 \end{gathered}
\end{equation}
In particular, \eqref{lambdaj} gives
\begin{equation}\label{height-contraction}
 \theta_{j+1}(c\rho_{j+1})^{2s}
 \le\widehat\gamma\lambda^{2s}\theta_j(c\rho_j)^{2s},
 \qquad \widehat\gamma\lambda^{2s}<1.
\end{equation}

\subsubsection{The tail estimate}
Split $\R^N\setminus K_{2\rho_j}$ into
$\R^N\setminus K_R$, $K_R\setminus K_{2\rho_0}$, and
$K_{2\rho_{i-1}}\setminus K_{2\rho_i}$, $1\le i\le j$.
The bound $|u|\le2M_{i-1}$ holds on each such annulus
throughout the time interval of $\widehat Q_j$. Therefore

\begin{align}\label{annular-decomposition}
 &\Tail(u;\widehat Q_j)\notag\\
 &=\int_{-\theta_j(c\rho_j)^{2s}}^{0}\!\int_{\R^N\setminus K_{2\rho_j}}\frac{|u(y,t)|}{|y|^{N+2s}}\,\dd y\,\dd t\notag\\
 &\le\int_{-\theta_j(c\rho_j)^{2s}}^{0}\!\int_{\R^N\setminus K_R}\frac{|u(y,t)|}{|y|^{N+2s}}\,\dd y\dd t+\int_{-\theta_j(c\rho_j)^{2s}}^{0}\!\int_{K_R\setminus K_{2\rho_0}}\frac{|u(y,t)|}{|y|^{N+2s}}\,\dd y\dd t\notag\\
&\qquad+\sum_{i=1}^{j}\int_{-\theta_j(c\rho_j)^{2s}}^{0}\!\int_{K_{2\rho_{i-1}}\setminus K_{2\rho_i}}\frac{|u(y,t)|}{|y|^{N+2s}}\,\dd y\dd t\notag\\
 &\le\int_{-\theta_j(c\rho_j)^{2s}}^0
       \int_{\R^N\setminus K_R}
       \frac{|u(y,t)|}{|y|^{N+2s}}\dd y\dd t+C\theta_j(c\rho_j)^{2s}
       \left(M_0\rho_0^{-2s}+\sum_{i=1}^jM_{i-1}\rho_i^{-2s}\right).
\end{align}
Using $M_{i-1}\le(1-\eta)^{-(j-i+1)}M_j$ and
$\rho_j/\rho_i=\lambda^{j-i}$, we have
\begin{align}
 \theta_j(c\rho_j)^{2s}M_0\rho_0^{-2s}
 &\le c^{2s}M_j^q
       \left(\frac{\lambda^{2s}}{1-\eta}\right)^j,
       \label{outer-annulus}\\
 \theta_j(c\rho_j)^{2s}\sum_{i=1}^jM_{i-1}\rho_i^{-2s}
 &\le\frac{c^{2s}M_j^q}{1-\eta}
       \sum_{i=1}^j\left(\frac{\lambda^{2s}}{1-\eta}\right)^{j-i}.
       \label{inner-annuli}
\end{align}
Impose $\lambda^{2s}/(1-\eta)\le1/2$ and sum the series.
Putting \eqref{outer-annulus},\eqref{inner-annuli} back in \eqref{annular-decomposition},
\begin{equation*}
 \Tail(u;\widehat Q_j)
 \le\int_{-\theta_j(c\rho_j)^{2s}}^0
       \int_{\R^N\setminus K_R}
       \frac{|u(y,t)|}{|y|^{N+2s}}\dd y\dd t
       +Cc^{2s}M_j^q.
\end{equation*}
Using $(a+b)^{1/q}\le C(q)(a^{1/q}+b^{1/q})$ and choosing
$$c\le\Bigl(\frac{\eta}{4\overline\gamma C}\Bigr)^{\frac{q}{2s}},$$  we obtain
from \eqref{recursive} that
\begin{align}\label{contracting-amplitude}
 2M_{j+1}
 \le\left(1-\frac\eta2\right)2M_j+C\left(\int_{-\theta_j(c\rho_j)^{2s}}^0
       \int_{\R^N\setminus K_R}
       \frac{|u(y,t)|}{|y|^{N+2s}}\dd y\dd t\right)^{1/q}.
\end{align}

\subsubsection{The discrete estimate}
The time intervals are nested by \eqref{height-contraction}.
Summing \eqref{contracting-amplitude} and using
\eqref{initial-amplitude}, we first obtain
\begin{equation}\label{amplitude-uniform}
 \begin{split}
 2M_j
 &\le2M_0+
 C\Tail\bigl(u;Q(R,\theta_0(c\rho_0)^{2s})\bigr)^{1/q}
 \le2\widetilde\gamma M_0,
 \end{split}
\end{equation}
where $\widetilde\gamma\ge1$ depends only on the $\data$.

For $j\ge2$, we now iterate \eqref{contracting-amplitude}
only from $\lfloor j/2\rfloor$ to $j$.
By \eqref{height-contraction}, for
$\lfloor j/2\rfloor\le i\le j$,
\[
 \theta_i(c\rho_i)^{2s}
 \le\theta_0(c\rho_0)^{2s}
       (\widehat\gamma\lambda^{2s})^{\lfloor j/2\rfloor}.
\]
The geometric sum of the coefficients of the tail terms
is bounded by
\[
 \sum_{i=\lfloor j/2\rfloor}^{j-1}
       \left(1-\frac\eta2\right)^{j-1-i}\le\frac2\eta.
\]
Consequently, \eqref{amplitude-uniform} gives
\begin{align}\label{discrete1}
 2M_j
 &\le2\widetilde\gamma M_0
       \left(1-\frac\eta2\right)^{j-\lfloor j/2\rfloor}
       \notag\\
 &\quad+C\left(
       \int_{-\theta_0(c\rho_0)^{2s}
                    (\widehat\gamma\lambda^{2s})^{\lfloor j/2\rfloor}}^0
       \int_{\R^N\setminus K_R}
       \frac{|u(y,t)|}{|y|^{N+2s}}\dd y\dd t
                   \right)^{1/q}.
\end{align}
This estimate holds at every completed step, including the
first stopping index. If the near-zero alternative persists,
both terms tend to zero as $j\to\infty$. The initial radius
$\rho_0$ remains fixed throughout the argument.

\subsection{Reduction of oscillation away from zero}
\label{sec:away-normalization}
Let $\ell$ be the first index at which \eqref{near-measure}
fails. Replacing $u$ with $-u$ if necessary, we may assume
\[
 |\{u\le M_\ell\}\cap Q_{c\rho_\ell}(\theta_\ell)|
       \le\nu|Q_{c\rho_\ell}(\theta_\ell)|.
\]
For $\ell\ge1$, \eqref{tailj} gives
\begin{align}\label{tailell}
 \Tail(u;\widehat Q_\ell)
 \le\left(\frac2{\overline\gamma}\right)^qM_\ell^q
       +\frac{Cc^{2s}}{1-\eta}M_\ell^q.
\end{align}
By \eqref{bar-gamma}, the first coefficient is at most
$1/(4\gamma_{**})$. Select $$c\leq \min\left\{\Bigl(\frac{\eta}{4\overline\gamma C}\Bigr)^{\frac{q}{2s}},\left(\frac{1}{8C\gamma_{**}}\right)^\frac{1}{2s}\right\},$$
so that the second
coefficient has the same bound. For $\ell=0$, use
\eqref{initial-amplitude} and \eqref{first-tail}.
Thus, in either case,
$\gamma_{**}\Tail(u;\widehat Q_\ell)\le M_\ell^q/2$.
Lemma~\ref{De01}, with radius $c\rho_\ell$ and reference
cylinder $\widehat Q_\ell$, now yields
\begin{equation*}
 \frac14M_\ell\le u\le2M_\ell
       \quad\text{a.e. in }Q_{c\rho_\ell/2}(\theta_\ell).
\end{equation*}
The upper bound follows from \eqref{induction-bounds}.

Introduce
\[
 \tilde x=x,\qquad
 \tilde t=\frac{t}{(M_\ell/4)^{q-1}},\qquad
 \bar v(\tilde x,\tilde t)
       =\frac4{M_\ell}u\bigl(\tilde x,(M_\ell/4)^{q-1}\tilde t\bigr).
\]
Then
\begin{equation*}
 1\le\bar v\le8
       \quad\text{in }Q_{c\rho_\ell/2}(4^{q-1}),
\end{equation*}
and $\bar v$ is a weak solution to
\begin{equation*}
 \partial_{\tilde t}\bar v^q+
                  \widetilde{\mathscr L}\bar v=0.
\end{equation*}
Here
$\widetilde K(\tilde x,\tilde y,\tilde t)
=K(\tilde x,\tilde y,(M_\ell/4)^{q-1}\tilde t)$
satisfies \eqref{kernel} with the same constants.

\subsection{The oscillation estimate after normalization}
Set
\[
 c_q=q\min\{1,8^{q-1}\},\qquad C_q=q\max\{1,8^{q-1}\}.
\]
For two numbers $a,k\in[1,8]$, by \eqref{g-comparison} and \eqref{power-difference}, we have 
\begin{equation}\label{quadratic-comparison}
 \begin{gathered}
 \frac{c_q}{2}(a-k)_\pm^2
       \le qg_\pm(a,k)\le\frac{C_q}{2}(a-k)_\pm^2,\\
 c_q(a-k)_\pm\le(a^q-k^q)_\pm\le C_q(a-k)_\pm.
 \end{gathered}
\end{equation}
Since $\bar v$ satisfies a similar equation as \eqref{equ01} in $Q_{c\rho_\ell/2}(4^{q-1})$, we may derive an analogous energy estimate in $K_r\times(\tau_1,\tau_2)\subset Q_{c\rho_\ell/2}(4^{q-1})$ as that of Proposition~\ref{energy} by using  \eqref{quadratic-comparison} and assuming that the time derivative of the moving level satisfies $k'\ge0$ in subsolution case, and $k'\le0$ in the supersolution case.

\begin{align}\label{quadratic-energy}
 &\frac{c_q}{2}\int_{K_r}w_\pm(x,\tau_2)^2\zeta(x,\tau_2)^2\dd x
 + c_q\int_{\tau_1}^{\tau_2}\int_{K_r}
                   |{k}'|\zeta^2w_\pm\dd x\dd\tilde t\notag\\
 &\quad+\frac12\int_{\tau_1}^{\tau_2}\int_{K_r}\int_{K_r}
      |w_\pm\zeta(x)-w_\pm\zeta(y)|^2
                \widetilde K\dd y\dd x\dd\tilde t\notag\\
 &\quad+\int_{\tau_1}^{\tau_2}\int_{K_r}\int_{K_r}
      w_\mp(y)w_\pm(x)\zeta(x)^2
                \widetilde K\dd y\dd x\dd\tilde t\notag\\
 &\le\frac{C_q}{2}\int_{K_r}w_\pm(x,\tau_1)^2\zeta(x,\tau_1)^2\dd x
 +\frac{C_q}{2}\int_{\tau_1}^{\tau_2}\int_{K_r}
                       w_\pm^2|\partial_{\tilde t}\zeta^2|
                              \dd x\dd\tilde t\notag\\
 &\quad+\int_{\tau_1}^{\tau_2}\int_{K_r}\int_{K_r}
      |\zeta(x)-\zeta(y)|^2w_\pm(x)w_\pm(y)
                   \widetilde K\dd y\dd x\dd\tilde t\notag\\
 &\quad+\int_{\tau_1}^{\tau_2}\int_{K_r}
                \int_{\R^N\setminus K_r}
       w_\pm(y)\zeta(x)^2w_\pm(x)\widetilde K
                               \dd y\dd x\dd\tilde t.
\end{align}

Time arguments in the integrands are $\tilde t$ when suppressed.
This is the counterpart of \cite[Proposition 2.1]{Liao-24b}
for $p=2$, from which Lemmas 3.1, 3.2 and 3.4 in \cite{Liao-24b} can be reproduced. We also note that the reduction $\mu^-=0$ made there for convenience is not needed. On the other hand, unlike in \cite{Liao-24b}, the coefficients $c_q/2,\, C_q/2$ of the first terms on the left and right sides are generally different. This will affect reproducing Lemma 3.3 of that paper.  To get it, we need to keep the exact time energy terms 
$\int_{K_r}g_\pm(u,{k})\zeta(x,\tau)^2\dd x\big|_{\tau_1}^{\tau_2}$, before using
\eqref{quadratic-comparison}, and estimate them as in \eqref{exact-time-energy01}--\eqref{exact-time-energy05}. Then the remaining terms in \eqref{quadratic-energy} can be handled by following the estimates in \cite[Lemma 3.3]{Liao-24b}.

With these preparation, we may therefore follow the oscillation and tail iteration
in \cite[Sections 4.2--4.4]{Liao-24b}.

To this end, set
\begin{align}\label{radiusell}
 \widehat R=2^{-[\frac{(q-1)_-}{s}+2]}c\rho_\ell.
\end{align}
Then, \eqref{tailell} and $|u|\le2M_\ell$
in $Q_\ell$ give
\begin{equation*}
 \int_{-\widehat R^{2s}}^0
       \int_{\R^N\setminus K_{c\rho_\ell/2}}
       \frac{|\bar v(\tilde y,\tilde t)|}
            {|\tilde y|^{N+2s}}\dd\tilde y\dd\tilde t
       \le C.
\end{equation*}
Here $C$ depends only on the $\data$, including $q$, and
is independent of $\ell$.
Indeed, 
\begin{align*}
    &\int_{-\widehat R^{2s}}^0
       \int_{\R^N\setminus K_{c\rho_\ell/2}}
       \frac{|\bar v(\tilde y,\tilde t)|}
            {|\tilde y|^{N+2s}}\dd\tilde y\dd\tilde t\\
&=\left(\frac{M_\ell}{4}\right)^{-q}\int_{-\widehat R^{2s} ( M_\ell/4)^{q-1}}^0
       \int_{\R^N\setminus K_{c\rho_\ell/2}}
       \frac{|u( y, t)|}
            {| y|^{N+2s}}\dd y\dd t\\
&=\left(\frac{M_\ell}{4}\right)^{-q}\int_{-\widehat R^{2s} (M_\ell/4 )^{q-1}}^0
       \int_{\R^N\setminus K_{2\rho_\ell}}
       \frac{|u( y, t)|}
            {| y|^{N+2s}}\dd y\dd t\\
&\qquad\quad+\left(\frac{M_\ell}{4}\right)^{-q}\int_{-\widehat R^{2s}\left(M_\ell/{4}\right)^{q-1}}^0
       \int_{K_{2\rho_\ell}\setminus K_{c\rho_\ell/2}}
       \frac{|u( y, t)|}
            {| y|^{N+2s}}\dd y\dd t\\
&\leq \left(\frac{M_\ell}{4}\right)^{-q}\Tail(u;\widehat Q_\ell)+ C\widehat R^{2s}(c\rho_\ell/2)^{-2s}\leq C.
\end{align*}
Here in the second-to-last inequality we used $\widehat R^{2s}\left(M_\ell/{4}\right)^{q-1}\leq \theta_\ell (c\rho_\ell)^{2s}$ to enlarge the time interval in the first integral and control it by the corresponding tail, and the fact that $|u|\le2M_\ell$
in $Q_\ell$ to estimate the second integral. To get the last bound, we used \eqref{tailell} and \eqref{radiusell}.

Thus,
\cite[Section 4.4]{Liao-24b} gives, for some $\beta\in(0,s)$  and
sufficiently small $r/\widehat R$,
\begin{align}\label{normalized-modulus}
 \essosc_{Q_r(1)}\bar v
 \le C\left(\frac r{\widehat R}\right)^\beta +C\int_{-C(\widehat Rr)^s}^0
       \int_{\R^N\setminus K_{c\rho_\ell/2}}
       \frac{|\bar v(\tilde y,\tilde t)|}
            {|\tilde y|^{N+2s}}\dd\tilde y\dd\tilde t .
\end{align}
The constants depend only on the $\data$.
The range of $r$ is chosen so that the time interval in the
tail is contained in $(-\widehat R^{2s},0)$.

\subsubsection{Rescaling to $u$}
For any measurable spatial region $A$ and admissible $H>0$,
the change of variables gives
\begin{align}\label{tail-scaling}
  \int_{-H}^0\int_A
       \frac{|\bar v(\tilde y,\tilde t)|}
            {|\tilde y|^{N+2s}}\dd\tilde y\dd\tilde t =
 \left(\frac{M_\ell}{4}\right)^{-q}
 \int_{-(M_\ell/4)^{q-1}H}^0\int_A
       \frac{|u(y,t)|}{|y|^{N+2s}}\dd y\dd t .
\end{align}
For almost every time in the interval of $\widehat Q_\ell$,
the annular decomposition used in
\eqref{annular-decomposition} gives
\begin{align*} 
 \int_{K_R\setminus K_{c\rho_\ell/2}}
       \frac{|u(y,t)|}{|y|^{N+2s}}\dd y \le C\left(M_\ell\rho_\ell^{-2s}
             +M_0\rho_0^{-2s}
             +\sum_{i=1}^\ell M_{i-1}\rho_i^{-2s}\right)
       \le CM_\ell\rho_\ell^{-2s}.
\end{align*}
Indeed, after division by $M_\ell\rho_\ell^{-2s}$,
the last two terms are bounded by
\[
 \left(\frac{\lambda^{2s}}{1-\eta}\right)^\ell
 +\frac1{1-\eta}\sum_{i=1}^\ell
       \left(\frac{\lambda^{2s}}{1-\eta}\right)^{\ell-i},
\]
which is uniformly bounded. For $\ell=0$ the sum is empty.
Consequently, whenever the corresponding time interval lies
in that of $\widehat Q_\ell$, \eqref{tail-scaling} yields
\begin{equation}\label{annulus-rescale}
 \begin{split}
 &\frac{M_\ell}{4}\int_{-H}^0
       \int_{K_R\setminus K_{c\rho_\ell/2}}
       \frac{|\bar v(\tilde y,\tilde t)|}
            {|\tilde y|^{N+2s}}\dd\tilde y\dd\tilde t
       \le CM_\ell H\rho_\ell^{-2s}.
 \end{split}
\end{equation}
Consequently,
\begin{align*}
&\int_{-C(\widehat Rr)^s}^0
       \int_{\R^N\setminus K_{c\rho_\ell/2}}
       \frac{|\bar v(\tilde y,\tilde t)|}
            {|\tilde y|^{N+2s}}\dd\tilde y\dd\tilde t\\
&=\left(\frac{M_\ell}{4}\right)^{-q}
 \int_{-(M_\ell/4)^{q-1}C(\widehat Rr)^s}^0\int_{\R^N\setminus K_{c\rho_\ell/2}}
       \frac{|u(y,t)|}{|y|^{N+2s}}\dd y\dd t\\
&=\left(\frac{M_\ell}{4}\right)^{-q}
 \int_{-(M_\ell/4)^{q-1}C(\widehat Rr)^s}^0\int_{\R^N\setminus K_{R}}
       \frac{|u(y,t)|}{|y|^{N+2s}}\dd y\dd t\\
&\qquad+\left(\frac{M_\ell}{4}\right)^{-q}
 \int_{-(M_\ell/4)^{q-1}C(\widehat Rr)^s}^0\int_{K_R\setminus K_{c\rho_\ell/2}}
       \frac{|u(y,t)|}{|y|^{N+2s}}\dd y\dd t\\
&\stackrel{\eqref{annulus-rescale}}{\leq}\left(\frac{M_\ell}{4}\right)^{-q}
 \int_{-(M_\ell/4)^{q-1}C(\widehat Rr)^s}^0\int_{\R^N\setminus K_{R}}
       \frac{|u(y,t)|}{|y|^{N+2s}}\dd y\dd t + C(r/\rho_\ell)^s
\end{align*}

Substituting above inequality into \eqref{normalized-modulus}, and rescaling $\bar v$ back to $u$ yield
\begin{align}\label{original-modulus}
 &\essosc_{K_r\times(-(M_\ell/4)^{q-1}r^{2s},0)}u\notag\\
 &\quad\le CM_\ell\left(\frac r{\rho_\ell}\right)^\beta
       +C\left(\frac{M_\ell}{4}\right)^{1-q}
       \int_{-C(M_\ell/4)^{q-1}(\rho_\ell r)^s}^0
       \int_{\R^N\setminus K_R}
       \frac{|u(y,t)|}{|y|^{N+2s}}\dd y\dd t .
\end{align}
Here $r/\rho_\ell$ is sufficiently small in terms of the
$\data$; in particular, the time interval in the integral
is contained in that of $\widehat Q_\ell$.
All constants are independent of the stopping index $\ell$.

\subsection{An explicit modulus of continuity}

Decrease the exponent $\beta$ in \eqref{original-modulus},
if necessary, so that
\begin{equation}\label{common-exponent}
 \begin{gathered}
 0<\beta<\frac s2,\\
 2\beta\log\bigl[\lambda(1-\eta)^{(q-1)_+/(2s)}\bigr]
       >\log(1-\eta/2),\\
 4\beta\log\bigl[\lambda(1-\eta)^{(q-1)_+/(2s)}\bigr]
       >\log(\widehat\gamma\lambda^{2s}).
 \end{gathered}
\end{equation}
Such a choice is possible since both expressions on the
right are negative. These inequalities also imply
$1-\beta(q-1)/(4s)>0$ when $q\ge1$.

If the iteration stops at $\ell$, we set $M_k=M_\ell$ for
$k>\ell$ when defining the remaining cylinders $Q_k$.
Then
\[
 (1-\eta)^kM_0\le M_k\le\tilde\gamma M_0,\qquad
 \theta_{k+1}\rho_{k+1}^{2s}
 \le\widehat\gamma\lambda^{2s}\theta_k\rho_k^{2s}
 <\theta_k\rho_k^{2s},
\]
and the cylinders remain nested. This continuation is
used only to select the index below; \eqref{discrete1}
will be used only at indices not exceeding $\ell$.

Fix $r\in(0,\rho_0)$. In what follows, we will choose $j$ by partitioning the
radii when $0<q<1$, and the cylinder heights when $q\ge1$.
The observation cylinder will be $Q_{\sigma r}(\theta_0)$,
where $\sigma\in(0,1)$ is a structural constant. Now, we take on the task.

\smallskip\noindent\textit{Case 1: $0<q<1$.}
Choose $j\ge0$ so that
\begin{equation}\label{radius-partition}
 \rho_{j+1}<r\le\rho_j.
\end{equation}
Since $M_j\le\widetilde\gamma M_0$ and $q-1<0$,
\[
 \theta_j=M_j^{q-1}\ge\widetilde\gamma^{q-1}\theta_0.
\]
Choose $\sigma^{2s}\le\widetilde\gamma^{q-1}$. The spatial
and temporal inclusions then follow from
\[
 \sigma r\le\rho_j,\qquad
 \theta_0(\sigma r)^{2s}\le\theta_j\rho_j^{2s}.
\]
Thus $Q_{\sigma r}(\theta_0)\subset Q_j$.
The other inequality in \eqref{radius-partition} gives
$r/\rho_0>\lambda^{j+1}$.

\smallskip\noindent\textit{Case 2: $q\ge1$.}
The heights $\theta_j\rho_j^{2s}$ decrease to zero.
Choose $j\ge0$ so that
\begin{equation}\label{height-partition}
 \theta_{j+1}\rho_{j+1}^{2s}
       <\theta_0r^{2s}\le\theta_j\rho_j^{2s}.
\end{equation}
The right-hand inequality and $M_j\le\widetilde\gamma M_0$
give
\[
 r\le\widetilde\gamma^{(q-1)/(2s)}\rho_j.
\]
Choose $\sigma\le\widetilde\gamma^{(1-q)/(2s)}$.
Then $\sigma r\le\rho_j$ and
$\theta_0(\sigma r)^{2s}\le\theta_j\rho_j^{2s}$, so again
$Q_{\sigma r}(\theta_0)\subset Q_j$.
The left-hand inequality, together with
$M_{j+1}\ge(1-\eta)^{j+1}M_0$, yields
\[
 \theta_0r^{2s}>
 \theta_0\rho_0^{2s}
       \bigl[(1-\eta)^{q-1}\lambda^{2s}\bigr]^{j+1}.
\]
For $q=1$, all $\theta_j=1$, and \eqref{height-partition}
is exactly the radius partition \eqref{radius-partition}.

In either case, we have
\begin{equation}\label{index-lower-bound}
 Q_{\sigma r}(\theta_0)\subset Q_j,\qquad
 \frac r{\rho_0}>
       \bigl[\lambda(1-\eta)^{(q-1)_+/(2s)}\bigr]^{j+1}.
\end{equation}
We may decrease $\sigma$ further below. Put
$i=\lfloor j/2\rfloor$, and first take $r/\rho_0$
sufficiently small that $i\ge2$. Now we consider two possibilities.

\smallskip\noindent\textit{The iteration reaches $i$.}
Suppose that there is no stopping index before $i$.
By nesting, $Q_{\sigma r}(\theta_0)\subset Q_i$.
Using \eqref{common-exponent}, \eqref{index-lower-bound},
and $i=\lfloor j/2\rfloor$, we obtain
\[
 \left(1-\frac\eta2\right)^{i-\lfloor i/2\rfloor}
 \le C\left(1-\frac\eta2\right)^{j/4}
 \le C\left(\frac{\sigma r}{\rho_0}\right)^{\beta/2}.
\]
Likewise, for sufficiently small $r/\rho_0$,
\[
 (c\rho_0)^{2s}
       (\widehat\gamma\lambda^{2s})^{\lfloor i/2\rfloor}
 \le C\rho_0^{2s}(\widehat\gamma\lambda^{2s})^{j/4}
 \le\rho_0^{2s}\left(\frac{\sigma r}{\rho_0}\right)^\beta.
\]
The last inequality uses the strict choice of $\beta$ in
\eqref{common-exponent}; $\sigma$ is a fixed structural
constant. Applying \eqref{discrete1} at index $i$ gives
\begin{align*}
 \essosc_{Q_{\sigma r}(\theta_0)}u
 &\le CM_0\left(\frac{\sigma r}{\rho_0}\right)^{\beta/2}
        +C\left(
       \int_{-\theta_0\rho_0^{2s}(\sigma r/\rho_0)^\beta}^0
       \int_{\R^N\setminus K_R}
       \frac{|u(y,t)|}{|y|^{N+2s}}\dd y\dd t
                   \right)^{1/q}.
\end{align*}

\smallskip\noindent\textit{The iteration stops before $i$.}
Let $\ell<i$ be the first stopping index, and set
\[
 \tilde r=\sigma r\max\left\{1,
        \left(\frac{4M_0}{M_\ell}\right)^{(q-1)/(2s)}\right\}.
\]
Then
\begin{equation}\label{common-cylinder}
 Q_{\sigma r}(\theta_0)\subset
 K_{\tilde r}\times
       (-(M_\ell/4)^{q-1}\tilde r^{2s},0),\qquad
 \tilde r\le C\sigma r(M_0/M_\ell)^{(q-1)_+/(2s)}.
\end{equation}
We use the two partitions once more to estimate
$\tilde r/\rho_\ell$.

For $0<q<1$, \eqref{radius-partition} implies
\begin{equation}\label{stopping-radius-singular}
 \frac{\tilde r}{\rho_\ell}
 \le C\sigma\lambda^{j-\ell}
 \le C\sigma\lambda^{j/2}
 \le C\left(\frac r{\rho_0}\right)^{1/2}.
\end{equation}
Here we used $\ell<i\le j/2$ and
$r>\rho_{j+1}$.

For $q\ge1$, the continuation gives
$M_j=M_{j+1}=M_\ell$. Thus \eqref{height-partition} reads
\begin{equation}\label{stopped-height-partition}
 \left(\frac{M_\ell}{M_0}\right)^{(q-1)/(2s)}
       \lambda^{j+1}
 <\frac r{\rho_0}
 \le
 \left(\frac{M_\ell}{M_0}\right)^{(q-1)/(2s)}
       \lambda^j.
\end{equation}
Its upper bound and \eqref{common-cylinder} yield
\begin{equation*}
 \frac{\tilde r}{\rho_\ell}
       \le C\sigma\lambda^{j-\ell}.
\end{equation*}
In both cases, we may fix $\sigma\leq c 2^{-[\frac{(q-1)_-}{s}+2]}$ for
\eqref{original-modulus} to apply at radius $\tilde r$.

For $0<q<1$, \eqref{stopping-radius-singular} and
$M_\ell\le\widetilde\gamma M_0$ directly give
\[
 M_\ell\left(\frac{\tilde r}{\rho_\ell}\right)^\beta
       \le CM_0\left(\frac r{\rho_0}\right)^{\beta/2}.
\]
For $q\ge1$, the lower bound in
\eqref{stopped-height-partition} gives
\begin{align*}
 M_\ell\left(\frac{\tilde r}{\rho_\ell}\right)^\beta
 &\le CM_\ell\lambda^{\beta(j-\ell)}\notag\\
 &\le CM_0\left(\frac r{\rho_0}\right)^{\beta/2}
       \left(\frac{M_\ell}{M_0}\right)^{1-\beta(q-1)/(4s)}
       \lambda^{\beta(j/2-\ell)}\notag\\
 &\le CM_0\left(\frac r{\rho_0}\right)^{\beta/2}.
\end{align*}
The last step uses $\ell<i\le j/2$,
$M_\ell\le\widetilde\gamma M_0$, and
$1-\beta(q-1)/(4s)>0$.

The tail interval in \eqref{original-modulus} is controlled
as follows. For $0<q<1$, \eqref{height-contraction} and
\eqref{stopping-radius-singular} give
\[
 C(M_\ell/4)^{q-1}(\rho_\ell\tilde r)^s
 \le C\theta_0\rho_0^{2s}
       (\widehat\gamma\lambda^{2s})^\ell
       \left(\frac r{\rho_0}\right)^{s/2}.
\]
For $q\ge1$, \eqref{common-cylinder} gives instead
\[
 C(M_\ell/4)^{q-1}(\rho_\ell\tilde r)^s
 \le C\theta_0\rho_0^{2s}
       \left(\frac{M_\ell}{M_0}\right)^{(q-1)/2}
       \left(\frac{\rho_\ell}{\rho_0}\right)^s
       \left(\frac r{\rho_0}\right)^s.
\]
Since $\beta<s/2$, both bounds imply
\begin{equation*}
 C(M_\ell/4)^{q-1}(\rho_\ell\tilde r)^s
       \le\theta_0\rho_0^{2s}
               \left(\frac{\sigma r}{\rho_0}\right)^\beta
\end{equation*}
for all sufficiently small $r/\rho_0$.

For $q\ge1$, compare $M_\ell$ with
\[
 \Tail\bigl(u;Q(R,\theta_0\rho_0^{2s}
                       (\sigma r/\rho_0)^\beta)\bigr)^{1/q}.
\]
If $M_\ell$ does not exceed this quantity, use
$\essosc_{Q_{\sigma r}(\theta_0)}u\le4M_\ell$.
Otherwise,
\[
 M_\ell^{1-q}
 \Tail\bigl(u;Q(R,\theta_0\rho_0^{2s}(\sigma r/\rho_0)^\beta)\bigr)
 \le
 \Tail\bigl(u;Q(R,\theta_0\rho_0^{2s}(\sigma r/\rho_0)^\beta)\bigr)^{1/q}.
\]
For $0<q<1$, the tail term in \eqref{original-modulus}
is bounded using $M_\ell^{1-q}\le C M_0^{1-q}$.
Together with the spatial estimates, this treats every
$\ell<i$, including $\ell=0$.

The normalized tail is at most one by
\eqref{initial-amplitude}. Consequently, both possibilities
give
\begin{align*}
 \essosc_{Q_{\sigma r}(\theta_0)}u
 &\le CM_0\left(\frac{\sigma r}{\rho_0}\right)^{\beta/2}
       \notag\\
 &\quad+CM_0\left[
       \frac1{M_0^q}
       \int_{-\theta_0\rho_0^{2s}(\sigma r/\rho_0)^\beta}^0
       \int_{\R^N\setminus K_R}
       \frac{|u(y,t)|}{|y|^{N+2s}}\dd y\dd t
                    \right]^{\min\{1,1/q\}}.
\end{align*}
Here the fixed spatial factors involving $\sigma$ have
been absorbed into $C$. Relabeling the observation radius
$\sigma r$ as $r$ proves \eqref{main-modulus} for
$r<\sigma\rho_0$, after translation to $(x_o,t_o)$.
For the remaining $\sigma\rho_0\leq r<\rho_0$, the estimate follows from
$\essosc_Q u\le4M_0$ after increasing $C$.
The initial radius $\rho_0$ has remained fixed.
This completes the proof.\hfill$\square$

\bigskip

\noindent\textbf{Declaration of the use of AI}:
The author used ChatGPT (OpenAI) to assist with language polishing, manuscript formatting and readability improvements. All AI-assisted content was independently verified by the author, and the use of AI did not affect the mathematical content or conclusions. The author takes full responsibility for the manuscript.

\end{document}